\documentclass[reqno]{amsart} 

\begin{document} 

\title[Heun equation II]
{The Heun equation and the Calogero-Moser-Sutherland system II:
perturbation and algebraic solution} 

\author[Kouichi Takemura]
{Kouichi Takemura} 

\address{Kouichi Takemura \hfill\break
Department of Mathematical Sciences,
Yokohama City University, 22-2 Seto, Kanazawa-ku, Yokohama
236-0027, Japan}
\email{takemura@yokohama-cu.ac.jp}


\keywords{Heun equation, Calogero-Moser-Sutherland system, \hfill\break\indent
Inozemtsev model, perturbation, Kato-Rellich theory, trigonometric limit, 
Heun function, \hfill\break\indent algebraic solution}

\begin{abstract}
 We apply a method of perturbation for the $BC_1$ Inozemtsev model 
 from the trigonometric model and show the holomorphy of perturbation.
 Consequently, the convergence of eigenvalues and eigenfuncions which 
 are expressed as formal power series is proved.
 We investigate also the relationship between $L^2$ space and some 
 finite dimensional space of elliptic functions.
\end{abstract}

\maketitle

\numberwithin{equation}{section}
\newtheorem{prop}{Proposition}[section]
\newtheorem{thm}[prop]{Theorem}
\newtheorem{lemma}[prop]{Lemma}
\newtheorem{cor}[prop]{Corollary}

\section{Introduction}

In this paper, we report some properties for eigenvalues and
eigenfuncions of the $BC_1$ Inozemtsev model. Consequently, we obtain results
 on Heun function.

The $BC_1$ Inozemtsev model is a one-particle model of quantum mechanics whose
 Hamiltonian is
\begin{equation}
H= -\frac{d^2}{dx^2} + \sum_{i=0}^3 l_i(l_i+1)\wp (x+\omega_i),
\label{Ino00}
\end{equation}
where $\wp (x)$ is the Weierstrass $\wp$-function with periods $(1, \tau)$,
$\omega _0=0$, $\omega _1=1/2$, $\omega_2=(1+\tau )/2$, $\omega_3=\tau /2$
are half-periods, and $l_i$, $(i=0,1,2,3)$ are coupling constants.
This model is sometimes called the $BC_1$ elliptic Inozemtsev model, because
the potential is described by use of elliptic functions.

There are two evidences which ensure the importance of the $BC_1$ Inozemtsev model.
The first one is equivalence to Heun equation, which will be explained in
section \ref{HI}.
The other one is that $BC_1$ quantum Inozemtsev model is a special $(N=1)$ case of
$BC_N$ Inozemtsev model \cite{Ino}, which is a generic integrable quantum system
with $B_N$ symmetry. In fact, classification of integrable quantum systems with
$B_N$ symmetry was done by Ochiai, Oshima, and Sekiguchi \cite{OOS}, and it was
shown that integrable quantum system with $B_N$ symmetry is $BC_N$ Inozemtsev model
or its degenerate one.
It is known that $BC_N$ Inozemtsev system contains the well-known
 Calogero-Moser-Sutherland system with $B_N$ symmetry as a special case.

In this paper, we try to obtain physical eigenfunctions and eigenvalues of the
$BC_1$ Inozemtsev model, and investigate their properties. Here a ``physical''
eigenfunction means that it is contained in an appropriate Hilbert space, which
is often a space of square-integrable ($L^2$) functions.
Note that roughly speaking the ``physical'' eigenfunction corresponds to the Heun
function of Heun equation.

Applying a method of perturbation is a possible approach to this problem, which
was done in \cite{Tak,KT} for the Calogero-Moser-Sutherland system of type $A_N$.
Now we explain this method shortly. Elliptic functions have a period $\tau$.
By a trigonometric limit $p =\exp(\pi \sqrt{-1} \tau) \to 0$,
the Hamiltonian of the $BC_1$ elliptic Inozemtsev model tends to the Hamiltonian
 of the $BC_1$ Calogero-Moser-Sutherland model, and it is known that eigenvalues
and eigenstates of the $BC_1$ Calogero-Moser-Sutherland model are obtained
explicitly by use of Jacobi polynomials.

Based on eigenstates for the case $p=0$, we can obtain eigenvalues and eigenstates
for the $BC_1$ elliptic Inozemtsev model $(p\neq 0)$ as formal power series in $p$.
This procedure is sometimes called an algorithm of perturbation
(see section \ref{sect:formalpert}).
Generally speaking, convergence of the formal power series obtained by perturbation
is not guaranteed a priori, but for the case of $BC_1$ elliptic Inozemtsev model,
the convergence radius of the formal power series in $p$ is shown to be non-zero
(see Corollary \ref{cor:conv}), and this perturbation is holomorphic.
As a result, real-holomorphy of eigenvalues in $ p$ and completeness of
eigenfunctions is proved. Note that a partial result was obtained in part I
\cite{Tak1} by applying Bethe Ansatz.

There is another method to investigate eigenvalues and eigenstates of the $BC_1$
Inozemtsev model.
If the coupling constants $l_0$, $l_1$, $l_2$, $l_3$ satisfy some equation, the
Hamiltonian $H$ (see (\ref{Ino00})) preserves a finite dimensional space of doubly
 periodic functions which is related to the quasi-exact solvability \cite{Tur,GKO}.
On a finite dimensional space, eigenvalues are calculated by solving the
characteristic equation, which is an algebraic equation, and eigenfunctions are
obtained by solving linear equations.
In this sense, eigenvalues on a finite dimensional space are ``algebraic'', and
eigenvalues and eigenfunctions on a finite dimensional space would be more explicit
than ones on an infinite dimensional Hilbert space.

In this paper, we also investigate relationship between Hilbert spaces ($L^2$ spaces)
and invariant spaces of doubly periodic functions with respect to the action of
the Hamiltonian $H$.
In some cases, a finite dimensional invariant space becomes a subspace of the
Hilbert space. Then it is shown under some assumption that the set of eigenvalues
on the finite dimensional invariant space coincides with the set of small
eigenvalues from the bottom on the Hilbert space.

This paper is organized as follows.
In section \ref{HI}, the relationship between the Heun equation and the $BC_1$
Inozemtsev system is clarified. Next we consider a trigonometric limit and review
that eigenstates for the trigonometric model are given by hypergeometric (Jacobi)
polynomials. We also explain how to apply an algorithm of perturbation in order to
obtain formal eigenvalues and eigenfunctions for the $BC_1$ elliptic
Inozemtsev model. In section \ref{sec:pertu}, we prove convergence of the
algorithm of perturbation  by applying Kato-Rellich theory. We also obtain several
results related to Kato-Rellich theory. Although holomorphy of the eigenfunctions
$\tilde{v}_m(x,p)$ in $p$ as elements of $L^2$ space is shown by Kato-Rellich
theory, convergence of the eigenfunctions for each $x$ is not assured immediately.
In section \ref{sec:holomo}, we show uniform convergence and holomorphy of the
eigenfunctions $\tilde{v}_m(x,p)$ for $x$ on compact sets.

In section \ref{sec:algfcn}, finite dimensional invariant subspaces of doubly
periodic functions are investigated and relationship to the Hilbert space ($L^2$
space) is discussed. In section \ref{sec:nonneg}, we focus on the case
$l_0, l_1, l_2, l_3 \in \mathbb{Z}_{\geq 0}$.
In section \ref{sec:example}, examples are presented to illustrate results of this
paper (especially section \ref{sec:algfcn}).
In section \ref{sec:com}, we give some comments.
In section \ref{sec:app}, some propositions are proved and definitions and properties of elliptic functions are provided.

We note that some results of this paper are generalized to the case of the
$BC_N$ Inozemtsev model (see \cite{Takq}).

\section{Heun equation, trigonometric limit and algorithm of perturbation} \label{HI}

\subsection{Heun equation and Inozemtsev system}
It is known that the Heun equation admits an expression in terms of elliptic
functions and this expression is closely related to the $BC_1$ Inozemtsev system
\cite{Ron,OOS,OS,Tak1}. In this subsection, we will explain this.

Let us recall the Hamiltonian of the $BC_1$ Inozemtsev model
\begin{equation}
H= -\frac{d^2}{dx^2} + \sum_{i=0}^3 l_i(l_i+1)\wp (x+\omega_i),
\label{Ino}
\end{equation}
where $\wp (x)$ is the Weierstrass $\wp$-function with periods
$(1, \tau)$, $\omega _0=0$, $\omega _1=1/2$, $\omega_2=(1+\tau )/2$,
$\omega_3=\tau /2$ are half-periods, and $l_i$ $(i=0,1,2,3)$ are coupling
constants. Assume that the imaginary part of $\tau $ is positive.
Set
\begin{gather*}
 e_i=\wp(\omega_i) \quad (i=1,2,3), \quad a=\frac{e_2-e_3}{e_1-e_3} , \\
 \tilde{\Phi}(w)=w^{\frac{l_0+1}{2}}(w-1)^{\frac{l_1+1}{2}}(aw-1)^{\frac{l_2+1}{2}}.
 \end{gather*}
Note that $a$ is nothing but the elliptic modular function $\lambda (\tau)$.
We change a variable by
\begin{equation}
w=\frac{e_1-e_3}{\wp (x)-e_3}.
\label{wxtrans}
\end{equation}
Then
\begin{equation} \label{Inotrans}
\begin{aligned}
 \tilde{\Phi}(w)^{-1} \circ H \circ \tilde{\Phi}(w)
&=-4(e_1-e_3)\Big\{ w(w-1)(aw-1) \big(\frac{d}{dw}\big) ^2   \\
&\quad +\frac{1}{2}\big( \frac{2l_0+3}{w}+  \frac{2l_1+3}{w-1}
+  \frac{a(2l_2+3)}{aw-1} \big) \frac{d}{dw} \Big\} +a\alpha \beta w +\tilde{q}
 \Big\},
\end{aligned}
\end{equation}
where
$\tilde{q}= \big( \frac{a+1}{3}\sum_{i=0}^3 l_i(l_i +1)-a(l_0+l_2+2)^2
-(l_0+l_1+2)^2 \big)$, $\alpha =\frac{l_0+l_1+l_2+l_3+4}{2}$,
$\beta =\frac{l_0+l_1+l_2-l_3+3}{2}$.

Let $f(x)$ be an eigenfunction of $H$ with an eigenvalue $E$, i.e.,
\begin{equation}
(H-E) f(x)= \Big( -\frac{d^2}{dx^2} + \sum_{i=0}^3 l_i(l_i+1)\wp (x+\omega_i)-E\Big)
f(x)=0.
\label{InoEF}
\end{equation}
 From (\ref{Inotrans}) and (\ref{InoEF}), we obtain
\begin{equation}
\Big(\big(\frac{d}{dw}\big) ^2  + \big( \frac{l_0+\frac{3}{2}}{w}
+\frac{l_1+\frac{3}{2}}{w-1}+\frac{l_2+\frac{3}{2}}{w-\frac{1}{a}}\big)
\frac{d}{dw} +\frac{\alpha \beta w -q}{w(w-1)(w-\frac{1}{a})} \Big)\tilde{f}(w)=0,
\label{Heun}
\end{equation}
where $\tilde{f}(\frac{e_1-e_3}{\wp (x)-e_3})
\tilde{\Phi}(\frac{e_1-e_3}{\wp (x)-e_3})=f(x)$ and
$q=-\frac{1}{4a}\big( \frac{E}{e_1-e_3}+\tilde{q} \big)$.
Note that the condition
\begin{equation}
l_0+\frac{3}{2}+l_1+\frac{3}{2}+l_2+\frac{3}{2}=\alpha +\beta +1
\label{HeunFR}
\end{equation}
is satisfied.

Equation (\ref{Heun}) with condition (\ref{HeunFR}) is called the Heun equation
\cite{Ron,SL}. It has four singular points $0$, $1$, $a^{-1}$,  $\infty$, all
the singular points are regular.
The following Riemann's $P$-symbol show the exponents.
\begin{equation*}
P\begin{pmatrix}
0 & 1 & a^{-1} & \infty & & \\
0 & 0 & 0 & \alpha & z & q \\
-l_0-\frac{1}{2} & -l_1-\frac{1}{2} & -l_2-\frac{1}{2} & \beta & &
\end{pmatrix}
\end{equation*}

Up to here, we have explained how to transform the equation of the Inozemtsev model
into the Heun equation.
Conversely, if a differential equation of second order with four regular singular
points on a Riemann sphere is given, we can transform it into equation (\ref{Heun})
with condition (\ref{HeunFR}) with suitable $l_i$ $(i=0,1,2,3)$ and $q$ by
changing a variable $w \to \frac{a'w+b'}{c'w+d'}$ and a transformation
$f \to w^{\alpha_1}(w-1)^{\alpha_2}(w-a^{-1})^{\alpha_3}f$.
It is known that if $a \neq 0,1 $ then there exists a solution $\tau $ to the
equation $a=\frac{e_2-e_3}{e_1-e_3}$ ($e_i$ $(i=1,2,3)$ depend on $\tau $).
Thus the parameter $\tau$ is determined.
The values $e_1,e_2,e_3$ and $E$ are determined by turn.
Hence we obtain a Hamiltonian of $BC_1$ Inozemtsev model with an eigenvalue $E$
starting from a differential equation of second order with four regular singular
points on a Riemann sphere.

\subsection{Trigonometric limit} \label{sect:trig}

In this section, we will consider a trigonometric limit
$(\tau \to \sqrt{-1} \infty)$.
We introduce a parameter $p=\exp (\pi \sqrt{-1} \tau )$, then $p\to 0$
as $\tau \to \sqrt{-1} \infty$. (Note that the parameter $p$ is different
from the one in \cite{Tak1}.)

If $p \to 0$, then $e_1 \to \frac{2}{3}\pi^2$,
$e_2 \to -\frac{1}{3}\pi^2$, $e_3 \to -\frac{1}{3}\pi^2$,
and $a \to 0$, and the relation between $x$ and $w$ (see (\ref{wxtrans}))
tends to $w=\sin ^2 \pi x$ as $p \to 0$.
Set
\begin{gather}
 H_T=  -\frac{d^2}{dx^2} + l_0(l_0+1)\frac{\pi ^2}{\sin ^2\pi x} +l_1(l_1+1)
 \frac{\pi ^2}{\cos ^2\pi x} ,\label{triIno} \\
 L_T = w(w-1)\Big\{ \frac{d^2}{dw^2} + \big( \frac{l_0+\frac{3}{2}}{w}
 + \frac{l_1+\frac{3}{2}}{w-1} \big)\frac{d}{dw} +\frac{(l_0+l_1+2)^2}{4w(w-1)}
 \Big\} \label{gauss}.
\end{gather}
Then $H \to H_T -\frac{\pi^2}{3}\sum_{i=0}^3 l_i (l_i +1)$ and
equation (\ref{Heun}) tends to
$(L_T -\frac{E}{\pi^2}-\frac{1}{3}\sum_{i=0}^3 l_i (l_i +1) )\tilde{f}(w)=0$
as $p \to 0$.
The operator $H_T$ is nothing but the Hamiltonian of the $BC_1$
trigonometric Calogero-Moser-Sutherland model, and the equation
$(L_T-C)\tilde{f}(w)=0$ ($C$ is a constant) is a Gauss hypergeometric equation.

Now we solve a spectral problem for $H_T$ by using hypergeometric functions.
We divide into four cases, the case $l_0>0$ and $l_1>0$, the case $l_0>0$ and
$l_1=0$, the case $l_0=0$ and $l_1>0$, and the case $l_0=l_1=0$.
For each case, we set up a Hilbert space $\mathbf{H}$, find a dense eigenbasis,
and obtain essential selfadjointness of the gauge-transformed trigonometric
Hamiltonian. The Hilbert space $\mathbf{H}$ plays an important role to show
holomorphy of perturbation in $p$, which will be discussed in
section \ref{sec:pertu}.
We note that the case $l_0=0$ and $l_1>0$ comes down to the case $l_0>0$ and
$l_1=0$ by setting $x \to x+\frac{1}{2}$.

\subsubsection{The case $l_0>0$ and $l_1>0$} \label{sssec1}
Set
\begin{equation*}
\Phi(x)=(\sin \pi x)^{l_0+1}(\cos \pi x)^{l_1+1}, \quad
\mathcal{H}_T=\Phi(x)^{-1} H_T \Phi(x),
\end{equation*}
then the gauge transformed Hamiltonian $\mathcal{H}_T$ is expressed as
\begin{equation}
\mathcal{H}_T= -\frac{d^2}{dx^2}-2\pi\Big( \frac{(l_0 +1)\cos \pi x}{\sin \pi x}
- \frac{(l_1 +1)\sin \pi x}{\cos \pi x} \Big) \frac{d}{dx}+(l_0+l_1+2)^2\pi^2.
\label{JacobiHam}
\end{equation}
By a change of variable $w=\sin ^2 \pi x$, we have
\begin{equation}\label{Ghyper}
\begin{aligned}
\: & \mathcal{H}_T -\pi^2(2m+l_0+l_1+2)^2 \\
 & = -4\pi ^2  \Big\{ w(1-w)\frac{d^2}{dw^2} \\
 & \quad +\Big(\frac{2l_0+3}{2}-((l_0+l_1+2)+1)w\Big)\frac{d}{dw}+m(m+l_0+l_1+2) \Big\}
\end{aligned}
\end{equation}
for each value $m$.
Hence the equation $\mathcal{H}_T -\pi^2(2m+l_0+l_1+2)^2$ is transformed into
a hypergeometric equation.
Set
\begin{equation}
\psi_m(x) = \psi^{(l_0,l_1)}_m(x) = \tilde{c}_m G_m
\Big( l_0+l_1+2, \frac{2l_0+3}{2}; \sin ^2 \pi x\Big) \quad (m \in \mathbb{Z}_{\geq 0}),
\label{Jacobipol}
\end{equation}
where the function $ G_m(\alpha, \beta; w)=_2 \! F_1(-m,\alpha+m;\beta;w)$ is
the Jacobi polynomial of degree $m$ and
$$
\tilde{c}_m=\sqrt{\frac{\pi (2m+l_0+l_1+2)\Gamma (m+l_0+l_1+2)
\Gamma (l_0+m+\frac{3}{2})}{m! \Gamma (m+l_1
+\frac{3}{2}) \Gamma (l_0+\frac{3}{2})^2}}
$$
is a constant for normalization. Then
\begin{equation}
\mathcal{H}_T \psi_m(x)=\pi^2(2m+l_0+l_1+2)^2\psi_m(x).
\label{eqn:Htpsim}
\end{equation}
We define the inner products
\begin{equation}
\langle f,g\rangle =\int_{0}^{1} dx\overline{f(x)} g(x), \quad
 \langle f,g\rangle _{\Phi} =\int_{0}^{1} dx\overline{f(x)} g(x) |\Phi(x)|^2.
\label{innerprod}
\end{equation}
Then $\langle \psi_m(x), \psi_{m'}(x) \rangle _{\Phi} =\delta_{m,m'}$.
Set
\begin{equation} \label{Hilb1}
\begin{aligned}
\mathbf{H} = \Big\{& f: \mathbb{R} \to \mathbb{C} : \mbox{measurable with }
 \int_{0}^{1} |f(x)| ^2 |\Phi (x)|^2 dx<+\infty, \\
 & f(x)=f(x+1), \; f(x)=f(-x) \mbox{ a.e. }x \Big\}
\end{aligned}
\end{equation}
and define an inner product on the Hilbert space $\mathbf{H}$ by $\langle \cdot , \cdot  \rangle _{\Phi}$.
Then the space spanned by functions $\{ \psi_m(x) | m\in \mathbb{Z}_{\geq 0} \}$ is dense in $\mathbf{H}$.
For $f(x),g(x) \in \mathbf{H} \cap C^{\infty}(\mathbb{R} )$, we have
\begin{equation} \label{HTsymm}
\begin{aligned}
\langle \mathcal{H}_T f(x),g(x)\rangle _{\Phi}
&= \langle \left(\mathcal{H}_T f (x)\right)\Phi(x) ,g(x)\Phi(x)\rangle
 = \langle H_T \left( f(x) \Phi(x)\right) ,g(x)\Phi(x)\rangle \\
&= \langle f(x) \Phi(x) ,H_T \left( g(x)\Phi(x)\right) \rangle
 =\langle f(x), \mathcal{H}_T g(x)\rangle _{\Phi}.
\end{aligned}
\end{equation}
It follows that the operator $\mathcal{H}_T$ is essentially selfadjoint on
the space $\mathbf{H}$.

\subsubsection{The case $l_0>0$ and $l_1=0$} \label{sssec2}
Set
\begin{equation*}
\Phi(x)=(\sin \pi x)^{l_0+1}, \quad \mathcal{H}_T=\Phi(x)^{-1} H_T \Phi(x),
\end{equation*}
then the gauge transformed Hamiltonian is expressed as
\begin{equation*}
\mathcal{H}_T= -\frac{d^2}{dx^2}-2\pi
\Big(\frac{(l_0 +1)\cos \pi x}{\sin \pi x}\Big) \frac{d}{dx}+(l_0+1)^2\pi^2.
\end{equation*}
By a change of variable $w=\sin ^2 \pi x$, we have a hypergeometric differential
equation,
\begin{equation} \label{GGhyper}
\begin{aligned}
\: &\mathcal{H}_T -\pi^2(2m+l_0+1)^2\\
&= -4\pi ^2 \Big\{ w(1-w)\frac{d^2}{dw^2}+\Big(\frac{2l_0+3}{2}-((l_0+1)+1)w\Big)
\frac{d}{dw}+m(m+l_0+1) \Big\}.
\end{aligned}
\end{equation}
for each $m$. Set
\begin{equation}
\psi^G _m(x)=  \tilde{c}^G_m C^{l_0+1}_m (\cos \pi x) \quad
(m \in \mathbb{Z}_{\geq 0}),
\label{Gegenpol}
\end{equation}
where the function
$C^{\nu}_m (z)=\frac{\Gamma (m+2\nu)}{m! \Gamma(2\nu)}\, _2 \! F_1(-m,m+2\nu ;
 \nu+\frac{1}{2} ;\frac{1-z}{2})$ is the Gegenbauer polynomial of degree $m$
and $\tilde{c}^G _m=\sqrt{\frac{2^{2l_0+1}(m+l_0+1)m!
\Gamma (l_0+1)^2}{\Gamma (m+2l_0+2)}}$.
Then
\[
\mathcal{H}_T \psi^G_m(x)=\pi^2(m+l_0+1)^2\psi^G_m(x),
\]
and $\langle \psi^G_m(x), \psi^G_{m'}(x) \rangle _{\Phi} =\delta_{m,m'}$,
where the inner product is defined as (\ref{innerprod}) for
$\Phi (x)=(\sin \pi x )^{l_0+1}$.

There are relations between Gegenbauer polynomials and Jacobi polynomials.
More precisely, $\psi^G _{2m}(x) = \psi^{(l_0,-1)} _{m}(x)$ and
$\psi^G _{2m+1}(x) = (\cos \pi x) \psi^{(l_0,0)} _{m}(x)$
$(m \in \mathbb{Z}_{\geq 0})$. Set
\begin{equation} \label{Hilbge}
\begin{gathered}
\begin{aligned}
\mathbf{H} = \Big\{ &f: \mathbb{R} \to \mathbb{C} : \mbox{measurable with}
\int_{0}^{1} |f(x)| ^2 |\Phi (x)|^2 dx<+\infty, \\
&\ f(x)=f(x+2), \; f(x)=f(-x) \mbox{ a.e. }x \Big\},
\end{aligned}\\
\mathbf{H}_+ = \{ f \in \mathbf{H} | f(x)=f(x+1) \mbox{ a.e. }x  \} ,\\
\mathbf{H}_- = \{ f \in \mathbf{H} | f(x)=-f(x+1) \mbox{ a.e. }x  \},
\end{gathered}
\end{equation}
and inner products on the Hilbert space $\mathbf{H}$ and its subspaces
$\mathbf{H}_+$, $\mathbf{H}_-$ are given by $\langle \cdot , \cdot  \rangle _{\Phi}$.
Then $\mathbf{H}_+ \perp \mathbf{H}_-$ and
$\mathbf{H}= \mathbf{H}_+ \oplus \mathbf{H}_-$.

The space spanned by $\{ \psi^G _m(x) | m\in \mathbb{Z}_{\geq 0} \}$ is
dense in $\mathbf{H}$. For $f(x),g(x) \in \mathbf{H} \cap C^{\infty}(\mathbb{R} )$,
we have $ \langle \mathcal{H}_T f(x),g(x)\rangle _{\Phi}=\langle f(x),
\mathcal{H}_T g(x)\rangle _{\Phi}$ similarly to (\ref{HTsymm}),
and it follows that the operator $\mathcal{H}_T$ is essentially selfadjoint on
the space $\mathbf{H}$.

Similar results hold for subspaces $\mathbf{H}_+ $ and $\mathbf{H}_- $.
In fact the space spanned by functions $\{ \psi_m(x) | m\in 2\mathbb{Z}_{\geq 0} \}$
(resp. $\{ \psi_m(x) | m\in 2\mathbb{Z}_{\geq 0} +1 \}$) is dense in $\mathbf{H}_+$
(resp. $\mathbf{H}_-$) and the operator $\mathcal{H}_T$ is essentially
selfadjoint on the space $\mathbf{H}_+$ (resp. $\mathbf{H}_-$).

\subsubsection{The case $l_0=0$ and $l_1>0$} \label{sssec21}
Although the case $l_0=0$ and $l_1>0$ comes down to the case $l_0>0$ and
$l_1=0$ by setting $x \to x+\frac{1}{2}$, we collect results for the $l_0=0$
and $l_1>0$ case for convenience.
Set
\begin{equation*}
\Phi(x)=(\cos \pi x)^{l_1+1}, \quad  \mathcal{H}_T=\Phi(x)^{-1} H_T \Phi(x),
\end{equation*}
then the gauge transformed Hamiltonian is expressed as
\begin{equation*}
\mathcal{H}_T= -\frac{d^2}{dx^2}-2\pi \Big(\frac{(l_1 +1)\sin \pi x}{\cos \pi x}
\Big) \frac{d}{dx}+(l_1+1)^2\pi^2.
\end{equation*}
Now we set
\begin{equation}
\psi^{G'} _m(x)=  \tilde{c}^{G'}_m C^{l_1+1}_m (\sin \pi x) \quad
(m \in \mathbb{Z}_{\geq 0}), \label{Gegenpol2}
\end{equation}
where $C^{\nu}_m (z)$ is the Gegenbauer polynomial appeared in (\ref{Gegenpol}) and \\
$\tilde{c}^{G'} _m=\sqrt{\frac{2^{2l_1+1}(m+l_1+1)m! \Gamma (l_1+1)^2}
{\Gamma (m+2l_1+2)}}$, then
\begin{align*}
\mathcal{H}_T \psi^{G'}_m(x)=\pi^2(m+l_1+1)^2\psi^{G'}_m(x),
\end{align*}
and $\langle \psi^{G'}_m(x), \psi^{G'}_{m'}(x) \rangle _{\Phi} =\delta_{m,m'}$, where the inner product is defined as (\ref{innerprod}).
Set
\begin{align*}
\mathbf{H} &= \Big\{ f: \mathbb{R} \to \mathbb{C} : \mbox{measurable wiht }
\int_{0}^{1} |f(x)| ^2 |\Phi (x)|^2 dx<+\infty, \\
&\quad f(x)=f(x+2), \; f(x)=f(-x+1) \mbox{ a.e. }x \Big\}, \\
 \mathbf{H}_+ &= \{ f \in \mathbf{H} | f(x)=f(x+1) \mbox{ a.e. }x  \} ,\\
\mathbf{H}_- &= \{ f \in \mathbf{H} | f(x)=-f(x+1) \mbox{ a.e. }x  \}.
\end{align*}
Here the inner product on the Hilbert space $\mathbf{H}$ is given by
$\langle \cdot , \cdot  \rangle _{\Phi}$.
Then $\mathbf{H}_+ \perp \mathbf{H}_-$ and $\mathbf{H}= \mathbf{H}_+ \oplus \mathbf{H}_-$.
The space spanned by functions $\{ \psi_m^{G'}(x) | m\in \mathbb{Z}_{\geq 0} \}$
is dense in $\mathbf{H}$, and the operator $\mathcal{H}_T$ is essentially
selfadjoint on the space $\mathbf{H}$.

\subsubsection{The case $l_0=0$ and $l_1=0$} \label{sssec3}
In this case, the trigonometric Hamiltonian is $H_T =-\frac{d^2}{dx^2}$.
Set $\Phi(x)=1$, $\mathcal{H}_T= H_T =-\frac{d^2}{dx^2}$,
$ \psi_m(x) = \sqrt{2} \cos m\pi x$,
$\varphi _m(x)=\sqrt{2} \sin m\pi x$ $(m\in \mathbb{Z}_{\geq 1})$, and
$ \psi_0(x) =1$.
Then
\begin{gather*}
\mathcal{H}_T \psi_m(x)=\pi^2m^2\psi_m(x), \quad (m\in \mathbb{Z} _{\geq 0}), \\
 \mathcal{H}_T \varphi _m(x)=\pi^2m^2\varphi_m(x), \quad
(m\in \mathbb{Z} _{>0}),  \\
 \langle \psi_m(x), \psi_{m'}(x) \rangle _{\Phi}
=\langle \varphi_m(x), \varphi_{m'}(x) \rangle _{\Phi} = \delta_{m,m'} , \\
 \langle \psi_m(x), \varphi_{m'}(x) \rangle _{\Phi} = 0 ,
\end{gather*}
where the inner product is defined by (\ref{innerprod}).
Set
\begin{align*}
\mathbf{H} &=\Big\{ f: \mathbb{R} \to \mathbb{C} : \mbox{measurable with}
\int_{0}^{1} |f(x)| ^2  dx<+\infty, \\
&\qquad  f(x)=f(x+2)  \mbox{ a.e. }x \Big\},\\
 \mathbf{H}_1 &= \{ f(x) \in \mathbf{H} : f(x)=f(x+1), \; f(x)=f(-x) \mbox{ a.e. } x\},  \\
 \mathbf{H}_2 &= \{ f(x) \in \mathbf{H} : f(x)=f(x+1), \; f(x)=-f(-x)\mbox{ a.e. } x\}, \\
 \mathbf{H}_3 &= \{ f(x) \in \mathbf{H} : f(x)=-f(x+1), \; f(x)=f(-x)\mbox{ a.e. } x\}, \\
 \mathbf{H}_4 &= \{ f(x) \in \mathbf{H} : f(x)=-f(x+1), \; f(x)=-f(-x)\mbox{ a.e. } x\}.
\end{align*}
Then the spaces $\mathbf{H}_i$ are pairwise orthogonal and
$\mathbf{H} = \oplus _{i=1}^4 \mathbf{H}_i$.

The space spanned by the functions $\{ \psi_m(x) | m\in \mathbb{Z}_{\geq 0} \}$ and
$\{ \varphi_m(x) | m\in \mathbb{Z}_{\geq 1} \}$ is dense in $\mathbf{H}$ and the
operator $\mathcal{H}_T$ is essentially selfadjoint on the space $\mathbf{H}$,
and also on subspaces $\mathbf{H}_i$, $(i=1,2,3,4)$.

\subsection{Perturbation on parameters $a$ and $p(=\exp (\pi \sqrt{-1} \tau ))$}
\label{sect:formalpert} $ $

As was explained in section \ref{sect:trig}, eigenvalues and eigenfunctions of the
Hamiltonian $H$ are obtained explicitly for the case $p=0$.
In this section, we apply a method of perturbation and have an algorithm for
obtaining eigenvalues and eigenfunctions as formal power series in $p$.

Since the functions $\wp(x+\omega_i)$ $(i=0,1,2,3)$ admit expansions (\ref{wpth})
and (\ref{wpth1}), the Hamiltonian $H$ (see (\ref{Ino})) admits the expansion
\begin{equation}
H=H_T +C_T+\sum_{k=1}^{\infty} V_k(x) p^k,
\label{Hamilp0}
\end{equation}
where $H_T$ is the Hamiltonian of the trigonometric model defined in (\ref{triIno}),  $V_k(x)$ $(k \in \mathbb{Z}_{\geq 1})$ are even periodic functions with period $1$, and $C_T= -\frac{\pi^2}{3}\sum_{i=0}^3 l_i(l_i+1)$ is a constant.
Note that for each $k$ the function $V_k(x)$ is expressed as a finite sum of
$\cos 2\pi n x$ $(n=0,\dots ,k)$.

First we consider the case ($l_0 >0$, $l_1 > 0$), ($l_0 >0$, $l_1 =0$) or
($l_0 = 0$, $l_1 >0$).
Set
\begin{equation*}
v_m = \begin{cases}
\psi _m(x)\Phi(x)  & (  l_0>0 , \; l_1>0) \\
\psi ^G _m(x)\Phi(x) & ( l_0>0 , \; l_1=0) \\
\psi ^{G'} _m(x)\Phi(x) & ( l_0=0 , \; l_1>0)
\end{cases}
\end{equation*}
for $m\in \mathbb{Z}_{\geq 0}$. Then $v_m $ $(m\in \mathbb{Z}_{\geq 0})$ is a normalized eigenvector of $H_T$ on $\mathbf{H}$.
Let $E_m$ be the eigenvalue of $H_T$ w.r.t. the eigenvector $v_m$, i.e.,
\begin{equation*}
E_m=\begin{cases}
(2m+l_0+l_1+2)^2 & (  l_0>0 , \; l_1>0) \\
(m+l_0+1)^2 & ( l_0>0 , \; l_1=0) \\
(m+l_1+1)^2 & ( l_0=0 , \; l_1>0).
\end{cases}
\end{equation*}
Then we have $E_m\neq E_{m'}$, if $m \neq m'$.

We will determine eigenvalues $E_m(p) = E_m+C_T+\sum_{k=1}^{\infty} E_{m}^{\{k\}}p^k$ and normalized eigenfunctions $v_m(p)= v_m+ \sum_{k=1}^{\infty} \sum_{m'} c_{m,m'}^{\{k\}}v_{m'}p^k$ of the operator $H_T +C_T+\sum_{k=1}^{\infty} V_k(x) p^k$ as formal power series in $p$.
In other words, we will find $E_m(p) $ and $v_m(p)$ that satisfy equations
\begin{gather}
 (H_T +C_T +\sum_{k=1}^{\infty} V_k(x) p^k )v_m(p) = E_m(p)v_m(p) ,\label{Hpertexp}\\
 \langle v_m(p) ,  v_m(p) \rangle =1 , \label{Hpertexp1}
\end{gather}
as formal power series in $p$.

First we calculate coefficients of $\sum_{m'} d_{m,m'}^{\{k\}}v_{m'}=V_k(x) v_{m}$.
Since $V_k(x)$ is a finite sum of $\cos 2n\pi x$ $(n=0,\dots ,k)$ and the
eigenvector $v_m$ is essentially a hypergeometric polynomial, coefficients
$d_{m,m'}^{\{k\}}$ are obtained by applying the Pieri formula repeatedly.
For each $m$ and $k$, $d_{m,m'}^{\{k\}}\neq 0$ for finitely many $m'$.

Now we compute $E_{m}^{\{k\}}$ and $ c_{m,m'}^{\{k\}}$ for $k \geq 1$.
Set $c_{m,m'}^{\{0\}} =\delta _{m,m'}$. By comparing coefficients of $v_{m'}p^k$,
it follows that conditions (\ref{Hpertexp}, \ref{Hpertexp1}) are equivalent to
following relations
\begin{align}
& c_{m,m'}^{\{ k \}}= \frac{\sum_{k'=1}^{k}
( \sum _{m''}c_{m,m''}^{\{ k-k' \} }d_{m'',m'}^{\{ k' \}})-\sum_{k'=1}^{k-1} c_{m,m'}^{\{ k-k' \} }E_{m}^{\{ k' \}} }
{E_m-E_{m'}} ,\quad (m' \neq m) \label{nondeg} \\
& c_{m,m}^{\{ k \}}=-\frac{1}{2} \Big(\sum_{k'=1}^{k-1} \sum_{m''}
c_{m,m''}^{\{ k' \}}c_{m,m''}^{\{ k-k' \} } \Big) ,\\
& E_m^{\{ k \}} = \sum_{k'=1}^{k}\sum_{m''}c_{m,m''}^{\{ k-k' \} }d_{m'',m}^{\{ k' \}}- \sum_{k'=1}^{k-1} c_{m,m}^{\{ k-k' \} } E_m^{\{ k' \}}.\label{nondeg2}
\end{align}
Note that the denominator of (\ref{nondeg}) is non--zero because of non-degeneracy
of eigenvalues.
Then numbers  $c_{m,m'}^{\{ k \}}$ and $E_m^{\{ k \}}$ are determined recursively
from (\ref{nondeg} - \ref{nondeg2}) uniquely.
It is shown recursively that for each $m$ and $k$,
$\# \{ m'' | \: c_{m,m''}^{\{ k \}} \neq 0 \}$ is finite and the sums
on (\ref{nondeg} - \ref{nondeg2}) in parameters $m''$ are indeed finite sums.

Therefore we obtain ``eigenvalues'' $E_m(p)$ and ``eigenfunctions'' $v_m(p)$ of
the operator $H$ as formal power series in $p$. At this stage, convergence is not
discussed.

Now consider the case $l_0=l_1=0$. Though there is degeneracy of eigenvalues
on the full Hilbert space $\mathbf{H}$, the degeneracy disappears when the action
of the Hamiltonian is restricted on $\mathbf{H}_i$ for each $i \in \{0,1,2,3 \}$
and the calculation of perturbation works compatibly on each space $\mathbf{H}_i$.
Hence the calculation is valid for the case $l_0=l_1=0$.

Let us discuss perturbation for equation (\ref{Heun}) on the parameter
$a\left( =\frac{e_2-e_3}{e_1-e_3} \right)$. We consider the case $l_0>0$
and $l_1>0$. Write
\begin{equation}
 L=(aw-1)L_T+ aw(w-1)\big(l_2+\frac{3}{2}\big)\frac{d}{dw}+aq_1w+\mathcal{E},
\label{eqnL}
\end{equation}
where
\begin{gather}
 L_T=  \left( w(w-1)\left(\frac{d^2}{dw^2}+\left( \frac{l_0+\frac{3}{2}}{w}+\frac{l_1+\frac{3}{2}}{w-1}\right) \frac{d}{dw}\right) +\frac{(l_0+l_1+2)^2}{4}\right), \label{eqnLT} \\
 q_1= \frac{(l_0+l_1+l_2-l_3+3)(l_0+l_1+l_2+l_3+4)-(l_0+l_1+2)^2}{4}, \nonumber \\
 \mathcal{E}=\frac{E}{4(e_1-e_3)}-\frac{(l_0+l_2+2)^2}{4}a-\sum_{i=0}^3\frac{l_i(l_i+1)e_3}{4(e_1-e_3)} .\nonumber
\end{gather}
Then the equation $L \tilde{f}(w)=0$ is equivalent to (\ref{Heun}).

We are going to find eigenvalues and eigenfunctions of $L$ as perturbation on $a$.
Set $\tilde{\psi}_m(w)={u}_m= _2\! F_1 (-m,m+l_0+l_1+2;\frac{2l_0 +3}{2};w)$, where
$\sin^2\pi x= w$.
We will use the following relations later which can be found in \cite{Ron}.
\begin{gather*}
 L_T \tilde{\psi}_m(w) = \frac{(2m+l_0+l_1+2)^2}{4} \tilde{\psi}_m(w), \\
 w \tilde{\psi}_m(w)= A_m \tilde{\psi}_{m+1}(w)+ B_m \tilde{\psi}_{m}(w)+C_m \tilde{\psi}_{m-1}(w), \\
 w(w-1) \frac{d}{dw} \tilde{\psi}_m(w)= A'_m \tilde{\psi}_{m+1}(w)+ B'_m \tilde{\psi}_{m}(w)+C'_m \tilde{\psi}_{m-1}(w).
\end{gather*}
where
\begin{gather*}
 A_m= -\frac{(m+l_0+l_1+2)(m+\frac{2l_0+3}{2})}{(2m+l_0+l_1+2)(2m+l_0+l_1+3)},\\
 A'_m= -\frac{m(m+l_0+l_1+2)(m+\frac{2l_0+3}{2})}{(2m+l_0+l_1+2)(2m+l_0+l_1+3)},\\
 B_m= \frac{2m(m+l_0+l_1+2)+\frac{2l_0+3}{2}(l_0+l_1+1)}{(2m+l_0+l_1+1)
 (2m+l_0+l_1+3)}, \\
 B'_m= \frac{m(m+l_0+l_1+2)(l_0-l_1)}{(2m+l_0+l_1+1)(2m+l_0+l_1+3)},  \\
 C_m= -\frac{m(m+\frac{2l_1+1}{2})}{(2m+l_0+l_1+2)(2m+l_0+l_1+1)}, \\
 C'_m= \frac{m(m+l_0+l_1+2)(m+\frac{2l_1+1}{2})}{(2m+l_0+l_1+2)(2m+l_0+l_1+1)}.
\end{gather*}
Set $\mathcal{E}_m= \frac{(2m+l_0+l_1+2)^2}{4}$,
 $\mathcal{E}_m(a)=\mathcal{E}_m+\sum_{k=1}^{\infty} 
 \mathcal{E}_{m}^{\{k\}} a^k$,\\
${u}_m(a)= {u}_m+ \sum_{k=1}^{\infty} \sum_{m'} \tilde{c}_{m,m'}^{\{k\}}{u}_{m'}a^k$,
and 
$\tilde{c}_{m,m}^{\{k\}}=0$ for $k\geq 1$ and all $m$. We will determine
$\mathcal{E}_{m}^{\{k\}}$ and $\tilde{c}_{m,m'}^{\{k\}}$ $(k \geq 1)$ to satisfy
$L \tilde{f}(w)=0$. Substituting in (\ref{eqnL}) with
$\mathcal{E}= \mathcal{E}_m (a)$, the following relations are shown:
\begin{gather}
 \tilde{c}_{m,m'}^{\{k\}} = \frac{1}{\mathcal{E}_m -\mathcal{E}_{m'}}
 \left( \tilde{c}_{m,m'+1}^{\{k-1\}} \tilde{C}_{m'+1}+ \tilde{c}_{m,m'}^{\{k-1\}}
 \tilde{B}_{m'}+ \tilde{c}_{m,m'-1}^{\{k-1\}} \tilde{A}_{m'-1} \right)
 \;\; (m' \neq m), \label{Tcjjp} \\
 \mathcal{E}_{m}^{\{k\}}= \tilde{c}_{m,m+1}^{\{k-1\}} \tilde{C}_{m+1}
 + \delta_{k-1,0} \tilde{B}_{m}+ \tilde{c}_{m,m-1}^{\{k-1\}} \tilde{A}_{m-1} , \label{TEjk}
\end{gather}
where
\begin{align*}
  \tilde{A}_{m}&= -(l_2+\frac{3}{2})A'_m-\mathcal{E}_mA_m-q_1A_m  \\
& = \frac{(m+l_0+l_1+2)(m+\frac{2l_0+3}{2})(m+\frac{l_0+l_1+l_2-l_3+3}{2})
(m+\frac{l_0+l_1+l_2+l_3+4}{2})}{(2m+l_0+l_1+2)(2m+l_0+l_1+3)} ,
\end{align*}
\begin{align*}
 \tilde{B}_{m}&= -(l_2+\frac{3}{2})B'_m-\mathcal{E}_mB_m-q_1B_m  \\
 &=-\Big(2m(m +l_0+l_1+2)+\frac{(2l_0+3)(l_0+l_1+1)}{2})
 \big(m(m+l_0+l_1+2)\\
 &\quad +\frac{(l_0+l_1+l_2-l_3+3)(l_0+l_1+l_2+l_3+4)}{4}\big)\Big)\\
 &\quad\div \Big((2m+l_0+l_1+3)(2m+l_0+l_1+1)\Big)\\
&\quad -\frac{(l_2+\frac{3}{2})(l_0-l_1)(m+l_0+l_1+2)m}{(2m+l_0+l_1+3)
(2m+l_0+l_1+1)} ,
\end{align*}
\begin{align*}
 \tilde{C}_{m}&= -(l_2+\frac{3}{2})C'_m-\mathcal{E}_mC_m-q_1C_m   \\
&=   \frac{m(m+\frac{2l_1+1}{2})(m+\frac{l_0+l_1-l_2+l_3+1}{2})
(m+\frac{l_0+l_1-l_2-l_3}{2})}{(2m+l_0+l_1+2)(2m+l_0+l_1+1)}.
\end{align*}
Solving the recursive equations (\ref{Tcjjp},\ref{TEjk}), $\mathcal{E}_m(a)$
and ${u}_m(a)$ are obtained as formal power series in $a$.
By expanding $e_1$, $e_2$, and $e_3$ as series in $a$, eigenvalues of the operator
$L$ are obtained as formal series in $a$.

Now we compare two expansions of eigenvalues
($E_m(p) \Leftrightarrow \mathcal{E}_m (a)$) and eigenfunctions
($v_m(p) \Leftrightarrow {u}_m(a)$).
 From the formula $a=16p\prod_{n=1}^{\infty}\big( \frac{1+p^{2n}}{1+p^{2n-1}} \big)^8$
(see \cite[\S 21.7]{WW}), it follows that $a$ is holomorphic in $p$ near $p=0$
 and admits an expansion $a=16p +O(p^2)$. Hence the formal power series ${u}_m(a)$
and $\mathcal{E}_m (a)$ are expressed as the formal power series in $p$, and
coefficients of $p^k$ on ${u}_m(a)$ (resp. $\mathcal{E}_m (a)$) are expressed as
linear combinations of coefficients of $a^l$ $(l\leq k)$ on ${u}_m(a)$
(resp. $\mathcal{E}_m (a)$).
Set $\widehat{v}_m(p)= (1-aw)^{-(l_2+1)/2}{u}_m(a)$ and
$\widehat{E}_m (p)=\mathcal{E}_m (a)$, then they also satisfy equation
(\ref{Hpertexp}) as formal power series in $p$.
Since the coefficients are determined by the recursive relations uniquely,
it follows that $\widehat{E}_m (p)= E_m(p)$ and $\widehat{v}_m(p) = C_m(p) v_m(p)$,
where $C_m(p)$ is a formal power series in $p$.
Note that $C_m(p)$ appears from a difference of the normalization.

In summary, the perturbation on the variable $p$ is equivalent to the one on the
variable $a$.

\section{Perturbation on the $L^2$ space} \label{sec:pertu}

Throughout this section, assume $l_0\geq 0$ and $l_1\geq 0$.

\subsection{Holomorphic perturbation}
In this subsection, we will use definitions and propositions written in
Kato's book \cite{Kat} freely. The main theorem in this subsection is
Theorem \ref{mainthmKato}. As an application, we show convergence of the formal
 power series of eigenvalues $E_m(p)$ (resp. $\mathcal{E}_{m}(a)$) in $p$
(resp. $a$) which are calculated by the algorithm of perturbation explained in
section \ref{sect:formalpert}.

We denote the gauge-transformed Hamiltonian of the $BC_1$ Inozemtsev model
by $H(p)$, i.e.,
\begin{equation}
H(p)= \Phi(x)^{-1} \circ \Big( - \frac{d^2}{dx^2} + \sum_{i=0}^3 l_i(l_i+1)\wp
(x+\omega_i) \Big) \circ \Phi(x),
\label{eqn:gtrans}
\end{equation}
where
$$
\Phi(x)= \begin{cases}
(\sin \pi x)^{l_0+1}(\cos \pi x)^{l_1+1} & (l_0>0, l_1>0); \\
(\sin \pi x)^{l_0+1} & (l_0>0, l_1=0 );\\
(\cos \pi x)^{l_1+1} & (l_0=0, l_1>0 );\\
1 & (l_0=0, l_1=0) ,
\end{cases}
$$
and $p=\exp(\pi \sqrt{-1} \tau)$.
Let $V_k(x)$ be the functions  in (\ref{Hamilp0}), $C_T$ be the constant in
(\ref{Hamilp0}), and $\mathcal{H}_T(=\Phi(x)^{-1} H_T \Phi(x))$ be the
gauge-transformed Hamiltonian of the $BC_1$ Calogero-Moser-Sutherland model.
Then the operator $H(p)$ is expanded as
\begin{equation*}
H(p)= \mathcal{H}_T+C_T+ \sum_{k=1}^{\infty} V_k(x)p^k.
\end{equation*}
Note that $H(0)=\mathcal{H}_T+C_T$. The functions $V_k(x)$ satisfy the following
lemma.

\begin{lemma} \label{lem:Vkx}
Let $s$ be a real number satisfying $s>1$. Then there exists a constant $A$
such that $V_k(x)\leq As^k$ for all $k \in \mathbb{Z}_{\geq 1}$ and $x \in \mathbb{R}$.
\end{lemma}

\begin{proof}
The functions $V_k(x)$ are determined by relations (\ref{wpth}).
Write
\begin{equation} \label{majser}
\begin{aligned}
 C(p)&=\sum_{k=1}^{\infty}V_kp^k\\
&=8\pi^2 \Big( \sum_{n=1}^{\infty} (l_0(l_0+1) +l_1(l_1+1))
\frac{2np ^{2n}}{1-p ^{2n}}\\
&\quad + \left( l_2(l_2+1)+l_3(l_3+1) \right) \big( \frac{n(p^{n}+p^{2n})}{1-p ^{2n}} \big)
 \Big) .
\end{aligned}
\end{equation}
Then the function $\sum _{k=1}^{\infty } V_k(x)$ is uniformly evaluated by
coefficients of majorant series $C(p)=\sum_{k=1}^{\infty}V_kp^k$ for $x\in \mathbb{R}$.
Since the convergence radius of (\ref{majser}) in $p$ is $1$, we obtain the lemma.
\end{proof}

\begin{prop}
The operator $H(p)$ $(-1<p<1)$ is essentially selfadjoint on the Hilbert space
$\mathbf{H}$.
\end{prop}

\begin{proof}
Fix the parameter $p$ $(-1<p<1)$.
Set $W(p)=H(p)-H(0)$. Then the operator $W(p)$ acts as multiplication of a
real-holomorphic periodic function of $x$ with period $1$. We denote the
 multiplication operator $W(p)$ by $W(p;x)$.

 From (\ref{majser}), an inequality $W(p;x)\leq C(p)$ is shown for all $x\in \mathbb{R}$.
Hence the boundedness $ \| W(p) f \| \leq C(p)\| f \|$ for $f\in \mathbf{H}$ is
shown. By definition of the inner product on $\mathbf{H}$, the operator $W(p)$
is shown to be symmetric on the whole space $\mathbf{H}$.
In sections \ref{sssec1}-\ref{sssec3}, essential selfadjointness of the
gauge-transformed trigonometric Hamiltonian $\mathcal{H}_T$ is illustrated.
Thus the operator $H(0)$ is also essentially selfadjoint.
Therefore essential selfadjointness of the operator $H(p)=H(0)+W(p)$ is derived
from the stability theorem for essential-selfadjoint operators
(see \cite[V-\S 4.1]{Kat}).
\end{proof}

Let $\tilde{T}$ denote the closure of a closable operator $T$.
Then $\tilde{H}(p)$ $(-1<p<1)$ is the unique extension of $H(p)$ to the
selfadjoint operator. Since the symmetric operator $W(p)$ is defined on the
whole space $\mathbf{H}$, the domain of the operator $\tilde{H}(p)$ $(-1<p<1)$
 coincides with the one of $\tilde{H}(0)$.

\begin{prop} \label{prop:holA}
The operators $\tilde{H}(p)$ form a holomorphic family of type (A) for $-1<p<1$.
\end{prop}

\begin{proof}
Let $s$ be a real number satisfying $s>1$. From Lemma \ref{lem:Vkx}, we have
$\| V_k(x) f \|\leq As^k \| f\|$, because the measure of an interval $(0,1)$ is $1$.
On the other hand, the operators $V_k(x)$ are defined on the whole space
$\mathbf{H}$ and symmetric.
Hence the proposition follows from Kato-Rellich theorem \cite[VII-\S 2.2,
Theorem 2.6]{Kat}.
\end{proof}

For the case $p=0$, the operator $\tilde{H}(0)$ coincides with the gauge-transformed
Hamiltonian of the trigonometric $BC_1$ Calogero-Moser-Sutherland model up to
constant. All eigenfunctions of the operator $\tilde{H}(0)$ were obtained
explicitly in section \ref{sect:trig}.
The spectrum $\sigma (\tilde{H}(0))$ contains only isolated spectra and the
multiplicity of each spectrum is $1$ or $2$.
Then the resolvent $R(\zeta, \tilde{H}(0))=(\zeta -\tilde{H}(0))^{-1}$ is compact
for $\zeta \not \in \sigma (\tilde{H}(0))$.

 From Theorem 2.4 in \cite[VII-\S 2.1]{Kat} and Proposition \ref{prop:holA} in this
 paper, we obtain the following statement.

\begin{prop} \label{prop:compresol}
The operator $\tilde{H}(p)$ has compact resolvent for all $p$ such that $-1<p<1$.
\end{prop}

Let $\sigma (\tilde{H}(p))$ be the spectrum of the operator $\tilde{H}(p)$.
 From Theorem 6.29 in \cite[III-\S 6.8]{Kat} and Proposition \ref{prop:compresol}
in this paper, it follows that

\begin{prop} \label{prop:discr}
The spectrum $\sigma (\tilde{H}(p))$ contains only point spectra and it is discrete.
The multiplicity of each eigenvalue is finite.
\end{prop}

Combining Theorem 3.9 in \cite[VII-\S 3.5]{Kat}, Propositions \ref{prop:holA}
and \ref{prop:compresol} in this paper and the selfadjointness of $\tilde{H}(p)$,
the following theorem is shown.

\begin{thm} \label{mainthmKato}
All eigenvalues of the Hamiltonian $\tilde{H}(p)$ $(-1<p<1)$ can be represented
as $\tilde{E}_m  (p)$ $(m\in \mathbb{Z}_{\geq 0})$, which is real-holomorphic in
$p \in (-1,1)$, and $\tilde{E}_m (0)$ coincides with the $(m+1)$st smallest
eigenvalue of the trigonometric Hamiltonian $H(0)$, which was obtained in
section \ref{sect:trig} explicitly.
The eigenfunction $\tilde{v}_m (p,x)$ of the eigenvalue $\tilde{E}_m (p)$ is
holomorphic in $p\in (-1,1)$ as an element in the Hilbert space $\mathbf{H}$,
and eigenvectors $\tilde{v}_m (p,x)$ form a complete orthonormal family on
$\mathbf{H}$.
\end{thm}

As an application of the theorem, convergence of the formal power series of
eigenvalues in the variable $p$ obtained by the algorithm of perturbation is shown.

\begin{cor}
Let $E_m(p)$ $(m\in \mathbb{Z}_{\geq 0})$ be the eigenvalue of the Hamiltonian $H(p)$
of the $BC_1$ Inozemtsev model obtained by the algorithm of perturbation
(see section \ref{sect:formalpert}), and $v_m(p)$ be the eigenvector of the
eigenvalue $E_m(p)$.
If $|p|$ is sufficiently small then the power series $E_m(p)$ converges, and the
power series $v_m(p)$ converges as an element in the Hilbert space.
\end{cor}

\begin{proof}
 From Theorem \ref{mainthmKato}, for each $m\in \mathbb{Z}_{\geq 0}$ and $p\in (-1,1)$,
there exists real-holomorphic eigenvalues $\tilde{E}_m (p)$ and normalized
eigenfunctions $\tilde{v}_m (p,x)\in L^2$ which converge to the trigonometric ones
as $p\to 0$. Since the eigenvalues $\tilde{E}_m (p)$ and the eigenfunctions
$\tilde{v}_m (p,x)$ are holomorphic in $p$ near $p=0$, there exists
$\epsilon_m \in \mathbb{R}_{>0}$ such that $\tilde{E}_m (p)$ and $\tilde{v}_m (p,x)$
are expanded as the series in $p$ and they converge on $|p|<\epsilon _m$.
The convergence of $\tilde{v}_m (p,x)$ is as an element in $\mathbf{H}$.
Eigenvalues $\tilde{E}_m (p)$ and eigenfunctions $\tilde{v}_m (p,x)$ satisfy the
following relations for $|p|<\epsilon _m$,
\begin{gather*}
 \Big( H_T +C_T +\sum_{k=1}^{\infty} V_k(x) p^{k} \Big)
  \Phi(x) \tilde{v}_m (p,x)= \tilde{E}_m (p)\Phi(x)\tilde{v}_m (p,x) ,\\
 \langle \Phi(x) \tilde{v}_m (p,x), \Phi(x)\tilde{v}_m (p,x)\rangle =1 .
\end{gather*}
These equations are the same as (\ref{Hpertexp}, \ref{Hpertexp1}). From uniqueness
of the coefficients obtained by perturbation, it is seen that
$\tilde{E}_m (p)= E_m(p)$ and $\Phi(x)\tilde{v}_m (p,x)= v_m(p)$.
Hence the convergence of $E_m(p)$ and $v_m(p)$ are shown.
\end{proof}

We also obtain the holomorphy of  perturbation on the variable $a$ from
Theorem \ref{mainthmKato}, because $p$ is holomorphic in $a$ near $a=0$
(see section \ref{sect:formalpert}).

\begin{cor} \label{cor:conv}
Let $\mathcal{E}_m(a)$ $(m\in \mathbb{Z}_{\geq 0})$ (see section \ref{sect:formalpert})
be the eigenvalue and ${u}_m(a)$ be the eigenvector of the Heun operator $L$
(see (\ref{eqnL})) obtained by the algorithm for perturbation from $L_T$
(see (\ref{eqnLT})).
If $|a|$ is sufficiently small then the power series $\mathcal{E}_m(a)$ converges,
and the power series ${u}_m(a)$ converges as an element in $L^2$ space.
\end{cor}

\subsection{Properties of the eigenvalue}
In this subsection, we will show some properties of eigenvalues.
First we discuss the multiplicity of eigenvalues. Under some assumptions, it is seen that eigenvalues never stick together.
We also introduce some inequalities for eigenvalues.

\begin{thm} \label{thm:multone}
Assume $l_0\geq 1/2$ or $l_1\geq 1/2$.
Let $\tilde{E}_m(p)$ $(m\in \mathbb{Z}_{\geq 0})$ be the eigenvalues of $\tilde{H}(p)$ defined in Theorem \ref{mainthmKato}. Then $\tilde{E}_m(p)\neq \tilde{E}_{m'}(p)$ $(m\neq m')$ for $-1<p<1$. In other words, eigenvalues never stick together.
\end{thm}

\begin{proof}
Assume $\tilde{E}_m(p_0)= \tilde{E}_{m'}(p_0)$ for some $m<m'$ and $-1<p_0<1$.
Let us consider solutions of a differential equation
$(\tilde{H}(p_0) -\tilde{E}_m(p_0))g(x)=0$. Exponents of solutions at $x=0$
(resp. $x=1/2$) are $l_0+1$ and $-l_0$ (resp. $l_1+1$ and $-l_1$).  From the
condition ($l_0\geq 1/2$ or $l_1\geq 1/2$), solutions $g(x)$ satisfying the $L^2$
condition $\int_0^1 |g(x)|^2dx$ form a one-dimentional space.
Therefore the multiplicity of the operator $\tilde{H}(p_0)$ on the Hilbert space
$\mathbf{H}$ with the eigenvalue $\tilde{E}_m(p_0)$ is one.

 From Proposition \ref{prop:discr}, the spectrum $\sigma(\tilde{H}(p_0))$
is discrete. Hence there exists $\epsilon_1 \in \mathbb{R} _{>0}$ such that
$\tilde{H}(p_0)$ has exactly one eigenvalue in the interval
$(\tilde{E}_m(p_0)-\epsilon _1,\tilde{E}_m(p_0)+\epsilon _1)$.
Since the operators $\tilde{H}(p)$ $(-1<p<1)$ form a holomorphic family, there
exists $\epsilon_2 \in \mathbb{R} _{>0}$ such that $\tilde{H}(p)$ has exactly one
eigenvalue in the interval
$(\tilde{E}_m(p_0)-\epsilon _1/2,\tilde{E}_m(p_0)+\epsilon _1/2)$ for $p$ such
that $p_0-\epsilon_2 <p< p_0+\epsilon_2$ (see \cite[V-\S 4.3, VII-\S 3.1]{Kat}).
Set $f(p)=\tilde{E}_m(p)- \tilde{E}_{m'}(p)$. Then the function $f(p)$ is
real-holomorphic in $p\in (-1,1)$, $f(p_0)=0$, and $f(0)\neq 0$.
 From the identity theorem for the real-holomorphic function $f(p)$,
it follows that the point $p=p_0$ is an isolated zero, i.e., there exists
$\epsilon \in \mathbb{R} _{>0}$ such that $\tilde{E}_m(p)\neq \tilde{E}_{m'}(p)$
for $0<|p-p_0|<\epsilon$.

Hence if $|p-p_0|(\neq 0)$ is sufficiently small then values
$\tilde{E}_m(p),\tilde{E}_{m'}(p)$ belong to the interval
$(\tilde{E}_m(p_0)-\epsilon _1/2,\tilde{E}_m(p_0)+\epsilon _1/2)$ and
$\tilde{E}_m(p) \neq \tilde{E}_{m'}(p) $. This shows that $\tilde{H}(p)$ has no
less than two eigenvalues in the interval
$(\tilde{E}_m(p_0)-\epsilon _1,\tilde{E}_m(p_0)+\epsilon _1)$, and it contradicts.
Therefore, we obtain the theorem.
\end{proof}

\begin{cor} \label{cor:ineq}
Assume $l_0\geq 1/2$ or $l_1\geq 1/2$, then $\tilde{E}_m(p)< \tilde{E}_{m'}(p)$
for $m<m'$ and $-1<p<1$.
\end{cor}

\begin{proof}
 From the labelling of the eigenvalues $\tilde{E}_m(p)$, it follows that
$\tilde{E}_m(0)< \tilde{E}_{m'}(0)$ for $m<m'$. It is seen in
Theorem \ref{thm:multone} that values $\tilde{E}_m(p)$ and $\tilde{E}_{m'}(p)$
never stick together for $-1<p<1$. Hence $\tilde{E}_m(p)< \tilde{E}_{m'}(p)$ for
$m<m'$ and $-1<p<1$.
\end{proof}

Now we introduce another inequality.

\begin{thm}
Assume $l_2=l_3=0$. Let $\tilde{E}_m(p)$ $(m\in \mathbb{Z}_{\geq 0})$ be the
eigenvalues of $\tilde{H}(p)$ defined in Theorem \ref{mainthmKato}.
For each $m$, the eigenvalue $\tilde{E}_m(p)$ is monotonely increasing for $0<p<1$
 and monotonely decreasing for $-1<p<0$.
\end{thm}

\begin{proof}
Set $\tilde{H}(p)= \tilde{H}(0)+\sum_{k=1}^{\infty} V_{2k}(x)p^{2k}$. From the
formulas (\ref{wpth}, \ref{wpth1}), it is seen that $V_{2k}(x)\geq 0$ for all
$k \in \mathbb{Z}_{>0}$. Hence
$\frac{d}{dp}\left( \sum_{k=1}^{\infty} V_{2k}(x)p^{2k}\right)\geq 0$ for
$0\leq p<1$ and $\frac{d}{dp}\left( \sum_{k=1}^{\infty} V_{2k}(x)p^{2k}\right)\leq 0$
for $-1<p\leq 0$.

Let $\tilde{v}_m (p,x)$ be the normalized eigenfunction of $\tilde{H}(p)$
labelled in Theorem \ref{mainthmKato}. By definition we have
$\tilde{H}(p)\tilde{v}_m (p,x) = \tilde{E}_m(p)\tilde{v}_m (p,x)$ and
\begin{equation}
\langle \tilde{H}(p)\tilde{v}_m (p,x) , \tilde{v}_m (p,x) \rangle _{\Phi}
= \tilde{E}_m(p).
\label{Hpeigen}
\end{equation}
By the way, we have
$\langle \frac{d}{dp} \tilde{v}_m (p,x), \tilde{v}_m (p,x) \rangle _{\Phi}=0$,
which is obtained by differentiating the equality
$\langle \tilde{v}_m (p,x), \tilde{v}_m (p,x) \rangle _{\Phi}=1$ in $p$.

Let us differentiate the equality (\ref{Hpeigen}). From the right hand side, we
have $\frac{d}{dp}\tilde{E}_m(p)$. From the left hand side, we obtain
\begin{align*}
& \Big\langle \Big( \frac{d}{dp}\sum_{k=1}^{\infty} V_{2k}(x)p^{2k} \Big)
\tilde{v}_m (p,x) , \tilde{v}_m (p,x) \Big\rangle _{\Phi}
+ 2 \big\langle \tilde{H}(p) \tilde{v}_m (p,x) , \frac{d}{dp} \tilde{v}_m (p,x)
 \big\rangle _{\Phi} \\
& = \int_0^1\Big( \frac{d}{dp} \sum_{k=1}^{\infty} V_{2k}(x)p^{2k} \Big)
|\tilde{v}_m (p,x)|^2|\Phi(x)|^2 dx +2 \tilde{E}_m(p)
\big\langle \tilde{v}_m (p,x), \frac{d}{dp} \tilde{v}_m (p,x) \big\rangle _{\Phi} \\
& =  \int_0^1 \Big( \frac{d}{dp} \sum_{k=1}^{\infty} V_{2k}(x)p^{2k} \Big)
|\tilde{v}_m (p,x)|^2|\Phi(x)|^2 dx.
\end{align*}
Here we used the selfadjointness of $\tilde{H}(p)$.
Thus
\begin{equation*}
\frac{d}{dp}\tilde{E}_m(p)
= \int_0^1 \Big( \frac{d}{dp}\sum_{k=1}^{\infty} V_{2k}(x)p^{2k} \Big)
|\tilde{v}_m (p,x)|^2 |\Phi(x)|^2 dx.
\end{equation*}
 If $0<p<1$ then the integrand is nonnegative for all $x \in [0,1]$. Hence $\frac{d}{dp}\tilde{E}_m(p)\geq 0$.
The inequality $\frac{d}{dp}\tilde{E}_m(p)\leq 0$ for $-1<p<0$ is shown similarly.
\end{proof}

\section{Perturbation on the space of holomorphic functions} \label{sec:holomo}

In the previous section, we obtain eigenfunctions $\tilde{v}_m(p,x)$ of the $BC_1$
dimensional Inozemtsev model as elements in the Hilbert space.
On the other hand, for the case $l_0, l_1\in \mathbb{Z}_{\geq 0}$ and $l_2=l_3=0$,
holomorophy of the function $\tilde{v}_m(p,x)$ in the variable $x$ for sufficintly
small $|p|$ is shown by applying the Bethe Ansatz method \cite{Tak1}.

In this section, we show holomorphy of the function $\tilde{v}_m(p,x)$ in the
variable $x$ for sufficintly small $|p|$ for the case $l_0, l_1\in \mathbb{R} _{\geq 0}$,
$l_2,l_3 \in \mathbb{R}$. The proof of holomorphy in $x$ is similar to the one for the
elliptic Calogero-Moser-Sutherland model of type $A_N$, which was done
in \cite{KT}. The paper \cite{KT} would be helpful to understand contents in this
section.

Let $v_m$ be a normalized eigenfunction of the trigonometric Hamiltonian
$\tilde{H}(0)$ and $E_m$ be the corresponding eigenvalue which were obtained
explicitly in sections \ref{sssec1}-\ref{sssec3}.

The action of the resolvent $(\tilde{H}(0)-\zeta )^{-1}$ on the Hilbert space
$\mathbf{H}$ is written as
\begin{equation*}
(\tilde{H}(0)-\zeta )^{-1}\sum_m  c_mv_m= \sum_m (E_m-\zeta)^{-1} c_mv_m,
\end{equation*}
for $\sum_m  c_mv_m \in \mathbf{H}$. Hence if $\zeta \not \in \sigma (\tilde{H}(0) )$,
 then the operator $(\tilde{H}(0)-\zeta )^{-1}$ is bounded.
Set $W(p)=\tilde{H}(p) -\tilde{H}(0)$. Then the operator $W(p)$ is bounded and
$\| W(p) \| \to 0$ as $p \to 0$. If  $\|(\tilde{H}(0)- \zeta )^{-1} W(p) \|<1$
then the resolvent $(\tilde{H}(p)-\zeta )^{-1}$ is expanded as
\begin{equation} \label{Neus}
\begin{aligned}
 (\tilde{H}(p)-\zeta )^{-1}
&= \bigl(1+(\tilde{H}(0)-\zeta )^{-1}W(p)\bigr)^{-1}(\tilde{H}(0)-\zeta )^{-1} \\
&= \sum_{j=0}^{\infty}\bigl(-(\tilde{H}(0)-\zeta )^{-1}W(p)\bigr)^j
(\tilde{H}(0)-\zeta )^{-1}.
\end{aligned}
\end{equation}
Let $\Gamma $ be a circle which does not bump the set $\sigma(\tilde{H}(0))$
and let $r=\mathop{\rm dist}(\Gamma,\sigma(\tilde{H}(0)))$.
Then there exists $p_0>0$ such that $\|W(p)\|<r$ and $\Gamma $ does not bump the
 set $\sigma(\tilde{H}(p))$ for all $|p|<p_0$.
Set
\begin{equation*}
 P_\Gamma(p):=-\frac{1}{2\pi \sqrt{-1}}\int_\Gamma(\tilde{H}(p)-\zeta )^{-1} d\zeta.
\end{equation*}
Then the operator $P_\Gamma(p)$ is a projection to the space of linear combinations
of eigenfunctions in the Hilbert space $\mathbf{H}$ whose eigenvalues are inside
the circle $\Gamma$.
Fix $E_m \in \sigma (\tilde{H}(0))$.
Since the set $\sigma (\tilde{H}(0))$ is discrete,
we can choose a circle $\Gamma _m$ and $p_m \in \mathbb{R} _{>0}$
such that $\Gamma _m$ contains only one element
$E_m$ of the set $\sigma (\tilde{H}(0))$ inside  $\Gamma _m$
and the projection $P_m(p)=P_{\Gamma_m}(p)$ satisfies
$\| P_m(p)-P_m(0) \|<1$ for $|p|<p_m$.

For the case ($l_0>0 ,\: l_1>0$), ($l_0>0 ,\: l_1=0$) or ($l_0=0 ,\: l_1>0$),
multiplicity of every eigenvalue of the operator $\tilde{H}(0)$ on the Hilbert
space $\mathbf{H}$ is one. Then a function $P_m(p)v_m$  for each $m$ is an
eigenvector of the operator $\tilde{H}(p)$ and admits an expression
\begin{equation}
P_m(p)v_m= c_{(p,m)} \tilde{v}_m(p,x), \label{pmvmcmvm}
\end{equation}
($c_{(p,m)}$ is a constant) for sufficiently small $|p|$, because the operator
$\tilde{H}(p)$ preserves the space $P_m(p) \mathbf{H}$ and $\tilde{v}_m(p,x)$
is an eigenvector of  $\tilde{H}(p)$ with the eigenvalue $\tilde{E}_m(p)$
(see Theorem \ref{mainthmKato}).

For the case $l_0=l_1=0$, the Hilbert space $\mathbf{H}$ is decomposed into
$\mathbf{H}_1\oplus \mathbf{H}_2\oplus \mathbf{H}_3\oplus \mathbf{H}_4$
(see section \ref{sssec3}). For each space, there is no degeneracy of eigenvalues
for the operator $\tilde{H}(0)$. Hence the expression (\ref{pmvmcmvm}) is also
valid.

As for the expansion of $P_m(p)v_m$ in $v_{m'}$ $({m'} \in \mathbb{Z}_{\geq 0})$,
the following proposition which is analogous to \cite[Proposition 5.11]{KT}
is shown.

\begin{prop}  \label{prop:pmvm}
Let $|p|<1$. Write $P_m(p)v_m = \sum_{m'} s_{m,{m'}} v_{m'}$. For each $m$ and
$C\in \mathbb{R} _{>1}$, there exists $C'\in \mathbb{R}_{>0}$ and $p_{\star }$ such that
coefficients $s_{m,{m'}}$ satisfy $|s_{m,{m'}}| \leq C'(C|p|)^{|m-{m'}|/2}$ for all
${m'} \in \mathbb{Z}_{\geq 0}$ and $p$ $(|p|<p_{\star })$.
\end{prop}

We prove this proposition in section \ref{sec:app}.
To obtain holomorphy of the function $P_m(p)v_m$, we need the following proposition,
which we will prove in section \ref{sec:app}.

\begin{prop} \label{prop:psival}
Let $\psi _m(x)(=v_m)$ be the $(m+1)$st normalized eigenfunction of the
trigonometric gauge-transformed Hamiltonian $\mathcal{H}_T$.
Let $f(x)=\sum_{m=0}^{\infty}c_m \psi _m(x)$ be a function satisfying
$|c_m|<A R^m$ $(\forall m \in \mathbb{Z} _{\geq 0})$ for some $A \in \mathbb{R}_{>0}$ and
$R \in (0,1)$.
If $r'$ satisfies $0<r'<\frac{1}{2\pi} \log \frac{1}{R}$, then the power series
$\sum_{m=0}^{\infty}c_m \psi _m(x)$ converges uniformly absolutely inside a
zone $-r' \leq \Im x\leq r'$, where $\Im x$ is an imaginary part of the complex
number $x$.
\end{prop}

\subsection*{Remark}
By the relation $w=\sin^2 \pi x$, a zone $| \Im x |<\frac{1}{2\pi} \log \frac{1}{R}$
in the variable $x$ corresponds to a domain inside an ellipse in the variable $w$
with fori at $w=0,1$ and a major axis $r=\frac{R+R^{-1}}{2}$.
\smallskip

Holomorphy of the power series $v_m (p)$ in variables $x$ and $p$ is shown by
applying Propositions \ref{prop:pmvm} and \ref{prop:psival}.

\begin{thm}
Let $v_m (p)$ be the eigenfunction of the Hamiltonian $\tilde{H}(p)$ obtained by
the algorithm of perturbation (see section \ref{sect:formalpert}).
Then the power series $v_m (p)$ is an analytic function in variables $x$ and $p$
for $|\Im x|$, $|p|$: sufficiently small.
More precisely for each $r \in \mathbb{R}_{>0}$ there exists $p_0 \in \mathbb{R}_{>0}$ such
that the series  $v_m(p)$ converges absolutely uniformly for
$(p,x) \in [-p_0 , p_0]\times B_r$, where $B_r= \{ x\in \mathbb{C} | | \Im x | \leq r\}$.
\end{thm}

\begin{proof}
We prove for the case $l_0>0$ and $l_1>0$.
Write $P_m(p)v_m = \sum_{m'} s_{m,{m'}} v_{m'}$ and apply
Proposition \ref{prop:pmvm} for the case $C=2$. Then there exists $C'\in \mathbb{R}_{>0}$
and $p_{\star }$ such that coefficients $s_{m,{m'}}$ satisfy
$|s_{m,{m'}}| \leq C'(2|p|)^{|m-{m'}|/2}$ for all $p$, ${m'}$ such that
$|p|<p_{\star }$.

Let $p_0 = \min (p_{\star }/2, \exp(-5\pi r)/2)$. Then we have
$r \leq \frac{2}{5\pi } \log \frac{1}{\sqrt{2p_0}}
 < \frac{1}{2\pi } \log \frac{1}{\sqrt{2p_0}}$ and there exists
 $A \in\mathbb{R} _{>0}$ such that $|s_{m,{m'}}| \leq A (\sqrt{2|p|})^{m'}$ for
$|p|\leq p_{0}(<p_{\star })$ and $m' \in \mathbb{Z}_{\geq 0}$.
By Proposition \ref{prop:psival}, the series $\sum_{m'} s_{m,{m'}} v_{m'}$
converges absolutely uniformly for $(p,x) \in [-p_0 , p_0]\times B_r$.
Since the function $v _m(p)$ constructed in section \ref{sect:formalpert}
coincides with a function $P_m(p)v_m$ up to constant and the function
$P_m(p)v_m = \sum_{m'} s_{m,{m'}} v_{m'}$ converges absolutely uniformly for
$(p,x) \in [-p_0 , p_0]\times B_r$, we obtain absolutely uniformly convergence
of $v _m(p)$.

For the case ($l_0>0$ and $l_1=0$) or ($l_0=0$ and $l_1>0$) or $(l_0=l_1=0)$,
we can prepare alternative propositions to Propositions \ref{prop:pmvm}
and \ref{prop:psival}, and can prove the theorem similarly.
\end{proof}

\section{Algebraic eigenfunctions} \label{sec:algfcn}

\subsection{Invariant subspaces of doubly periodic functions} \label{sec:invsp}
If the coupling constants $l_0$, $l_1$, $l_2$, $l_3$ satisfy some equation,
the Hamiltonian $H$ (see (\ref{Ino})) preserves a finite dimensional space of
doubly periodic functions. In this section, we look into a condition for existence
of the finite dimensional invariant space of doubly periodic functions with respect
to the action of the Hamiltonian $H$ (see Proposition \ref{findim}). After that,
we investigate relationship between invariant spaces of doubly periodic functions
and $L^2$ spaces.

\begin{prop} \label{findim}
Let $\tilde{\alpha}_i$ be a number such that $\tilde{\alpha}_i= -l_i$ or
$\tilde{\alpha}_i= l_i+1$ for each $i\in \{ 0,1,2,3\} $.
Assume $-\sum_{i=0}^3 \tilde{\alpha}_i /2\in \mathbb{Z}_{\geq 0}$ and set
$d=-\sum_{i=0}^3 \tilde{\alpha}_i /2$.
Let $V_{\tilde{\alpha}_0, \tilde{\alpha}_1, \tilde{\alpha}_2, \tilde{\alpha}_3}$ be
the $d+1$ dimensional space spanned by
$\Big\{ (\frac{\sigma_1(x)}{\sigma(x)})^{\tilde{\alpha }_1}
(\frac{\sigma_2(x)}{\sigma(x)})^{\tilde{\alpha }_2}(\frac{\sigma_3(x)}{\sigma(x)}
)^{\tilde{\alpha }_3}\wp(x)^n \: | \:n=0, \dots ,d \Big\}$,
where $\sigma (x)$ is the Weierstrass sigma-function defined in (\ref{wpsigzeta})
and $\sigma_i(x)$ ($i=1,2,3$) are the co-sigma functions defined in (\ref{cosigma}).
Then the Hamiltonian $H$ (see (\ref{Ino})) preserves the space
$V_{\tilde{\alpha}_0, \tilde{\alpha}_1, \tilde{\alpha}_2, \tilde{\alpha}_3}$.
\end{prop}

\begin{proof}
Set $z=\wp(x)$, $\widehat{\Phi}(z)=(z-e_1)^{\tilde{\alpha}_1/2}
(z-e_2)^{\tilde{\alpha}_2/2}(z-e_3)^{\tilde{\alpha}_3/2}$, and
$\widehat{H}= \widehat{\Phi}(z)^{-1} \circ H \circ\widehat{\Phi}(z)$. Then
\begin{align*}
 \widehat{H}= & -4(z-e_1)(z-e_2)(z-e_3)
 \Big\{ \frac{d^2}{dz^2}+ \sum_{i=1}^3\frac{\tilde{\alpha}_i+\frac{1}{2}}{z-e_i}
 \frac{d}{dz}\Big\} \\
& +\big\{ -\left(\tilde{\alpha}_1+\tilde{\alpha}_2+\tilde{\alpha}_3-l_0\right)
\left(\tilde{\alpha}_1+\tilde{\alpha}_2+\tilde{\alpha}_3+l_0+1\right) z\nonumber \\
& + e_1(\tilde{\alpha}_2+\tilde{\alpha}_3)^2+e_2(\tilde{\alpha}_1
+\tilde{\alpha}_3)^2+e_3(\tilde{\alpha}_1+\tilde{\alpha}_2)^2 \big\} \nonumber .
\end{align*}
Let $\widehat{V}_{d+1}$ be the space of polynomials in $z$ with degree at most $d$.
From formulas (\ref{sigmai}), it is enough to show that the operator $\widehat{H}$
preserves the space $\widehat{V}_{d+1}$.

The action of the Hamiltonian is written as
\begin{equation}\label{ze2}
\begin{aligned}
 \widehat{H}(z-e_2)^r
&= -4\Big( (r+\gamma_1)(r+\gamma_2) (z-e_2)^{r+1}  \\
& \quad  +((e_2-e_3)(r+\tilde{\alpha}_2+\tilde{\alpha}_1)r+(e_2-e_1)
(r+\tilde{\alpha}_2+\tilde{\alpha}_3)r+q')(z-e_2)^r  \\
& \quad+ r\big(r+\tilde{\alpha}_2-\frac{1}{2}\big)(e_2-e_3)(e_2-e_1)(z-e_2)^{r-1}\Big),
\end{aligned}
\end{equation}
where $q'=-\frac{1}{4}\left(e_1(\tilde{\alpha}_2+\tilde{\alpha}_3)^2
+e_2(\tilde{\alpha}_1+\tilde{\alpha}_3)^2+e_3(\tilde{\alpha}_1
+\tilde{\alpha}_2)^2\right)+e_2\gamma_1\gamma_2$, \\
$\gamma_1=(\tilde{\alpha}_1+\tilde{\alpha}_2+\tilde{\alpha}_3-l_0)/2$ and
$\gamma_2=(\tilde{\alpha}_1+\tilde{\alpha}_2+\tilde{\alpha}_3+l_0+1)/2$.
Hence the operator $\widehat{H}$ preserves the space of polynomials in $z$.
Since $\tilde{\alpha}_0= -l_0$ or $l_0+1$, it follows that
$(r+\gamma_1)(r+\gamma_2)=0$ for $r=d$, $(d=-\sum_{i=0}^3 \tilde{\alpha}_i/2)$.
Hence the coefficient of $(z-e_2)^{r+1}$ on the right hand side of (\ref{ze2})
is zero for the case $r=d$.
Therefore the operator $\widehat{H}$ preserves the space $\widehat{V}_{d+1}$.
\end{proof}

\begin{prop} \label{prop:fintrig}
With the notation in Proposition \ref{findim},
assume $d=-\sum_{i=0}^3 \tilde{\alpha}_i/2 \in \mathbb{Z}_{\geq 0}$. For the
trigonometric case $(p=0)$, the eigenvalues of the Hamiltonian $H$ on the finite
 dimensional space
$V_{\tilde{\alpha}_0, \tilde{\alpha}_1, \tilde{\alpha}_2, \tilde{\alpha}_3}$
are written as $\{ \pi^2 (2r+ \tilde{\alpha}_0+ \tilde{\alpha}_1)^2-\frac{\pi^2}{3}
 \sum_{i=0}^3 l_i(l_i+1)\} _{r=0, \dots ,d}$.
\end{prop}

\begin{proof}
As $p\to 0$, we have $e_1\to 2\pi^2/3$, $e_2\to -\pi^2/3$, and
$e_3\to -\pi^2/3$.
In this case the coefficient of $(z-e_2)^{r-1}$ on the right hand side of (\ref{ze2}) is zero for all $r$. Hence the operator $H$ acts triangularly.
Then the eigenvalues appear on diagonal elements. By a straightforward calculation,
the coefficient of $(z-e_2)^{r}$ on the right hand side of (\ref{ze2})
is written as
\begin{equation*}
\pi^2 (2r+ \tilde{\alpha}_2+ \tilde{\alpha}_3)^2-\frac{\pi^2}{3}
 \sum_{i=0}^3 l_i(l_i+1).
\end{equation*}
Here we used relations $ \tilde{\alpha}_i^2- \tilde{\alpha}_i= l_i(l_i+1)$
$(i=0,1,2,3)$. Therefore the eigenvalues of the Hamiltonian $H$ on the finite
dimensional space
$V_{\tilde{\alpha}_0, \tilde{\alpha}_1, \tilde{\alpha}_2, \tilde{\alpha}_3}$ are
$\{ \pi^2 (2r+ \tilde{\alpha}_2+ \tilde{\alpha}_3)^2-\frac{\pi^2}{3}
\sum_{i=0}^3 l_i(l_i+1)\} _{r=0, \dots ,d}$.
By replacing $r \to d-r$, we obtain the proposition.
\end{proof}

\subsection{Algebraic eigenfunctions on the Hilbert space} \label{sec:algL2}

In this subsection we investigate a relationship between the Hilbert space
$\mathbf{H}$ and the finite dimensional space
$V_{\tilde{\alpha}_0, \tilde{\alpha}_1, \tilde{\alpha}_2, \tilde{\alpha}_3}$ which
was defined in the previous subsection.
Throughout this subsection, assume $l_0\geq 0$ and $l_1 \geq 0$.

Let us consider the case $\tilde{\alpha}_i \in \{ -l_i , l_i+1 \}$ $(i=0,1,2,3)$ and
$d=-\sum_{i=0}^3 \tilde{\alpha}_i/2 \in \mathbb{Z}_{\geq 0}$.
It is easily confirmed that if $\tilde{\alpha}_0 \geq 0$ and
$\tilde{\alpha}_1 \geq 0$ then any function $f(x)$ in
 $V_{\tilde{\alpha}_0, \tilde{\alpha}_1, \tilde{\alpha}_2, \tilde{\alpha}_3}$ 
is square-integrable on the interval $(0,1)$.

Set $\tilde{V}_{\tilde{\alpha}_0, \tilde{\alpha}_1, \tilde{\alpha}_2,
\tilde{\alpha}_3}= \{ \frac{f(x)}{\Phi(x)} | f(x) \in V_{\tilde{\alpha}_0,
\tilde{\alpha}_1, \tilde{\alpha}_2, \tilde{\alpha}_3 }\}$, where $\Phi(x)$
is the ground state of the trigonometric model which was defined in
section \ref{sect:trig}.

\begin{prop} \label{prop:finL2}
Let $\tilde{\alpha}_i \in \{ -l_i, l_i+1\}$ for $i=0,1,2,3$.
If  $\tilde{\alpha}_0 \geq 0$, $\tilde{\alpha}_1 \geq 0$ and
$d=-\sum_{i=0}^3 \tilde{\alpha}_i /2\in \mathbb{Z}_{\geq 0}$, then
$\tilde{V}_{\tilde{\alpha}_0, \tilde{\alpha}_1, \tilde{\alpha}_2, \tilde{\alpha}_3}
\subset \mathbf{H}$.
For the case ($l_0>0$ and $l_1=0$) or ($l_0=0$ and $l_1>0$) we have
$\tilde{V}_{\tilde{\alpha}_0, \tilde{\alpha}_1, \tilde{\alpha}_2, \tilde{\alpha}_3}
\subset \mathbf{H}_i$ for $i=+$ or $-$, and for the case $l_0=l_1=0$ we have
$\tilde{V}_{\tilde{\alpha}_0, \tilde{\alpha}_1, \tilde{\alpha}_2, \tilde{\alpha}_3}
\subset \mathbf{H}_i$ for $i=1,2,3$ or $4$.
\end{prop}

\begin{proof}
We consider four cases, i.e., the case $l_0>0$, $l_1>0$, the case $l_0>0$,
$l_1=0$ case, the case $l_0=0,\; l_1>0$, and the case $l_0=l_1=0$.

First, we prove the proposition for the case $l_0>0$, $l_1>0$.
The numbers $\tilde{\alpha}_0$ and $\tilde{\alpha}_1$ must be chosen as
$\tilde{\alpha}_0= l_0+1$ and $\tilde{\alpha}_1= l_1+1$ from the condition
 $\tilde{\alpha}_0 ,\tilde{\alpha}_1 ,l_0 ,l_1 >0$.
Now we check that, if $f(x) \in \tilde{V}_{\tilde{\alpha}_0, \tilde{\alpha}_1,
 \tilde{\alpha}_2, \tilde{\alpha}_3}$, then $f(x)$ satisfies the definition of the
Hilbert space $\mathbf{H}$ (see (\ref{Hilb1})). Square-integrability of the
function $f(x)$ follows from the condition $\tilde{\alpha}_0\geq 0$ and
$\tilde{\alpha}_1\geq 0$. Periodicity and symmetry of $f(x)$ follow from the
condition $-\sum_{i=0}^3 \tilde{\alpha}_i /2 \in \mathbb{Z} $. Hence $f(x)\in \mathbf{H}$.

For the case $l_0>0$ and $l_1=0$,  $\tilde{\alpha}_0$ and $\tilde{\alpha}_1$ are
chosen as $\tilde{\alpha}_0= l_0+1$ and ($\tilde{\alpha}_1= 0$ or $1$).
Then $\tilde{V}_{\tilde{\alpha}_0, \tilde{\alpha}_1, \tilde{\alpha}_2,
\tilde{\alpha}_3} \subset \mathbf{H}_i$ for $i=+$ or $-$ is shown similarly.
 Note that the sign of $i$ is determined by whether functions in
$\tilde{V}_{\tilde{\alpha}_0, \tilde{\alpha}_1, \tilde{\alpha}_2, \tilde{\alpha}_3}$
are periodic or antiperiodic (for details see section \ref{sec:L2l0l12}).

For the other cases, the proofs are similar.
\end{proof}

\subsection*{Remark}
A function in $\tilde{V}_{\tilde{\alpha}_0, \tilde{\alpha}_1, \tilde{\alpha}_2,
 \tilde{\alpha}_3}$ is multi-valued in general, and so we should specify branches
 of the function. In our case, the analytic continuation of the function near the
real line should be performed along paths passing through the upper half plane.
\smallskip

In the finite dimensional space $\tilde{V}_{\tilde{\alpha}_0, \tilde{\alpha}_1,
\tilde{\alpha}_2, \tilde{\alpha}_3}$, eigenvalues are calculated by solving a
characteristic equation, which is an algebraic equation, and eigenfunctions are
obtained by solving linear equations.
In this sense, eigenvalues in the finite dimensional space are ``algebraic''.

Now we figure out properties of the spaces
$\tilde{V}_{\tilde{\alpha}_0, \tilde{\alpha}_1, \tilde{\alpha}_2, \tilde{\alpha}_3}$
contained in the Hilbert space $\mathbf{H}$.
We divide into four cases.

\subsubsection{The case $l_0>0$ and $l_1>0$}
Let $\tilde{\alpha}_2 \in \{ -l_2, l_2+1 \}$ and
$\tilde{\alpha}_3 \in \{ -l_3, l_3+1 \}$.
 From Proposition \ref{prop:finL2}, if $d=-(l_0+l_1)/2-1-(\tilde{\alpha}_2
 +\tilde{\alpha}_3 )/2\in \mathbb{Z}_{\geq 0}$ then the $d+1$ dimensional vector space
$\tilde{V}_{l_0+1, l_1+1, \tilde{\alpha}_2, \tilde{\alpha}_3}$ is a subspace of the
Hilbert space $\mathbf{H}$. We will show that the set of eigenvalues of the
gauge-transformed Hamiltonian $H(p)$ (see (\ref{eqn:gtrans})) on the space
$\tilde{V}_{l_0+1, l_1+1, \tilde{\alpha}_2, \tilde{\alpha}_3}$ is the set of
small eigenvalues of $H(p)$ on the Hilbert space $\mathbf{H}$ from the bottom.

\begin{lemma}
For the trigonometric case $p=\exp(\pi \sqrt{-1} \tau) = 0$, if
$d=-(l_0+l_1)/2-1-\tilde{\alpha}_2 -\tilde{\alpha}_3 \in \mathbb{Z}_{\geq 0}$ then the
set of eigenvalues of the gauge-transformed Hamiltonian $H(0)$ on the finite
dimensional space
$\tilde{V}_{l_0+1, l_1+1,  \tilde{\alpha}_2, \tilde{\alpha}_3}|_{p=0}$ coincides
with the set of small eigenvalues of $H(0)$ on the Hilbert space $\mathbf{H}$ from
the bottom. In other words, the $m$-th smallest eigenvalue of $H(0)$ on $\mathbf{H}$
is also an eigenvalue on
$\tilde{V}_{l_0+1, l_1+1,  \tilde{\alpha}_2, \tilde{\alpha}_3}|_{p=0}$, if
and only if $m\leq d+1$.
\end{lemma}

\begin{proof}
 From Proposition \ref{prop:fintrig}, eigenvalues of the trigonometric Hamiltonian
$H$ on the finite dimensional space
$V_{\tilde{\alpha}_0, \tilde{\alpha}_1, \tilde{\alpha}_2, \tilde{\alpha}_3}$ are
written as $\{ \pi^2 (2r+l_0+l_1+2)^2-\frac{\pi^2}{3} \sum_{i=0}^3 l_i(l_i+1)
\} _{r=0, \dots ,d}$. On the other hand, from equality (\ref{eqn:Htpsim}) and the
limit $H \to H_T-\frac{\pi^2}{3} \sum_{i=0}^3 l_i(l_i+1)$ as $p\to 0$,
eigenvalues of the gauge-transformed trigonometric Hamiltonian
$H(0)=\Phi(x)H \Phi(x)^{-1}|_{p=0}$ on the Hilbert space $\mathbf{H}$ are written
as $\{ \pi^2 (2r+l_0+l_1+2)^2-\frac{\pi^2}{3} \sum_{i=0}^3 l_i(l_i+1)\}
 _{r\in \mathbb{Z}_{\geq 0}}$.
Therefore the lemma follows.
\end{proof}

 From the previous lemma, eigenvalues of the trigonometric gauge-transformed
Hamiltonian $H(0)$ on the finite dimensional space
$ V_{l_0+1, l_1+1,  \tilde{\alpha}_2, \tilde{\alpha}_3}$ are
$\{ \tilde{E}_r(0) \: | \: r=0,\dots ,d \}$, where the values $\tilde{E}_r(p)$ are
defined in Theorem \ref{mainthmKato}.
It is obvious that the eigenvalues of the operator $H(p)$ on the finite dimensional
space $\tilde{V}_{l_0+1, l_1+1,  \tilde{\alpha}_2, \tilde{\alpha}_3}$ are
continuous in $p$. Hence if $-1<p <1$ then the eigenvalues of the gauge-transformed
Hamiltonian $H(p)$ on the finite dimensional space
$\tilde{V}_{l_0+1, l_1+1,  \tilde{\alpha}_2, \tilde{\alpha}_3}$ coincide with
$\{ \tilde{E}_r(p)\: | \: r=0,\dots ,d \}$.
By applying Corollary \ref{cor:ineq} for the case $l_0\geq 1/2$ or the case
$l_1\geq 1/2$, we obtain the following theorem.

\begin{thm} \label{thm:lowest}
Assume $d=-(l_0+l_1)/2-1-(\tilde{\alpha}_2 +\tilde{\alpha}_3)/2 \in \mathbb{Z}_{\geq 0}$.
If ($l_0\geq 1/2$ and $l_1>0$) or ($l_0>0$ and $l_1\geq 1/2$), then the set of
 eigenvalues of the gauge-transformed Hamiltonian $H(p)$ on the finite dimensional
space $\tilde{V}_{l_0+1, l_1+1,  \tilde{\alpha}_2, \tilde{\alpha}_3}$  coincides
with the set of small eigenvalues of $H(p)$ on the Hilbert space $\mathbf{H}$ from
the bottom.
In other words, the $m$-th smallest eigenvalue of $H(p)$ on $\mathbf{H}$
is also an eigenvalue on $\tilde{V}_{l_0+1, l_1+1,  \tilde{\alpha}_2,
\tilde{\alpha}_3}$, if and only if $m\leq d+1$.
\end{thm}

\subsubsection{The case $l_0>0$ and $l_1=0$} \label{sec:L2l0l12}
Let $\tilde{\alpha}_2 \in \{ -l_2, l_2+1 \}$ and
$\tilde{\alpha}_3 \in \{ -l_3, l_3+1 \}$.
If $d=-(l_0+1)/2-(\tilde{\alpha}_2 +\tilde{\alpha}_3)/2 \in \mathbb{Z}_{\geq 0}$,
then the $d+1$ dimensional space $\tilde{V}_{l_0+1, 0, \tilde{\alpha}_2,
\tilde{\alpha}_3}$ is a subspace of the space $\mathbf{H}_+$, and if
$d=-(l_0+2)/2-(\tilde{\alpha}_2 +\tilde{\alpha}_3)/2 \in \mathbb{Z}_{\geq 0}$,
then the $d+1$ dimensional space $\tilde{V}_{l_0+1, 1, \tilde{\alpha}_2,
\tilde{\alpha}_3}$ is a subspace of the space $\mathbf{H}_-$.

Similarly to Theorem \ref{thm:lowest}, the following statement are shown.\\
$\bullet$ If $d=-(l_0+1)/2-(\tilde{\alpha}_2 +\tilde{\alpha}_3)/2 \in
\mathbb{Z}_{\geq 0}$, $l_0\geq 1/2$, and $-1<p<1$, then the set of eigenvalues of
the gauge-transformed Hamiltonian $H(p)$ (see (\ref{eqn:gtrans})) on the finite
dimensional space $\tilde{V}_{l_0+1, 0, \tilde{\alpha}_2, \tilde{\alpha}_3}$ is
the set of small eigenvalues of $H(p)$ on the space $\mathbf{H}_+$ from the bottom.\\
$\bullet $ If $d=-(l_0+2)/2-(\tilde{\alpha}_2 +\tilde{\alpha}_3)/2 \in
\mathbb{Z}_{\geq 0}$, $l_0\geq 1/2$, and $-1<p<1$, then the set of eigenvalues of
$H(p)$ on the finite dimensional space $\tilde{V}_{l_0+1, 1, \tilde{\alpha}_2,
\tilde{\alpha}_3}$ is the set of small eigenvalues of $H(p)$ on the Hilbert space
$\mathbf{H}_-$ from the bottom.

\subsubsection{The case $l_0=0$ and $l_1>0$}  \label{sec:L2l0l13}
Let $\tilde{\alpha}_2 \in \{ -l_2, l_2+1 \}$ and
$\tilde{\alpha}_3 \in \{ -l_3, l_3+1 \}$.
If $d=-(l_1+1)/2-(\tilde{\alpha}_2 +\tilde{\alpha}_3)/2 \in \mathbb{Z}_{\geq 0}$,
then $\tilde{V}_{0,l_1+1, \tilde{\alpha}_2, \tilde{\alpha}_3} \subset \mathbf{H}_+$
and $\dim \tilde{V}_{0,l_1+1, \tilde{\alpha}_2, \tilde{\alpha}_3} =d+1$.
If $d=-(l_1+2)/2-(\tilde{\alpha}_2 +\tilde{\alpha}_3)/2 \in \mathbb{Z}_{\geq 0}$, then
$\tilde{V}_{1,l_1+1, \tilde{\alpha}_2, \tilde{\alpha}_3} \subset \mathbf{H}_-$ and
$\dim \tilde{V}_{1,l_1+1, \tilde{\alpha}_2, \tilde{\alpha}_3} =d+1$.
We can confirm similar statements to the ones in section \ref{sec:L2l0l12}.

\subsubsection{The case $l_0=0$ and $l_1=0$}
Let $\tilde{\alpha}_2 \in \{ -l_2, l_2+1 \}$ and $\tilde{\alpha}_3 \in
 \{ -l_3, l_3+1 \}$.
If $d=-(\tilde{\alpha}_2 +\tilde{\alpha}_3)/2 \in \mathbb{Z}_{\geq 0}$, then
$\tilde{V}_{0, 0, \tilde{\alpha}_2, \tilde{\alpha}_3}\subset \mathbf{H}_1$,
$\tilde{V}_{1, 1, \tilde{\alpha}_2, \tilde{\alpha}_3}\subset \mathbf{H}_2$,
$\dim \tilde{V}_{0, 0, \tilde{\alpha}_2, \tilde{\alpha}_3} =d+1$, and
$\dim \tilde{V}_{1, 1, \tilde{\alpha}_2, \tilde{\alpha}_3}=d$.
If $d=-1/2-(\tilde{\alpha}_2 +\tilde{\alpha}_3)/2 \in \mathbb{Z}_{\geq 0}$, then
$\tilde{V}_{1, 0, \tilde{\alpha}_2, \tilde{\alpha}_3}\subset \mathbf{H}_4$,
$\tilde{V}_{0, 1, \tilde{\alpha}_2, \tilde{\alpha}_3}\subset \mathbf{H}_3$, and
$\dim \tilde{V}_{1, 0, \tilde{\alpha}_2, \tilde{\alpha}_3}
=\dim \tilde{V}_{0, 1, \tilde{\alpha}_2, \tilde{\alpha}_3}=d+1$.

\subsection{The case of nonnegative integral coupling constants}
\label{sec:nonneg}
If the coupling constants $l_0, l_1, l_2, l_3$ are nonnegative integers, the model
satisfies some special properties. Specifically, it admits the Bethe Ansatz method
\cite{Tak1} and the potential has the finite-gap property \cite{TV,Smi,Tak3}.
In this subsection, we reproduce several results
more explicitly for the case $l_0, l_1, l_2, l_3$ are nonnegative integers.
Note that some results were obtained in \cite{Tak1}.

Throughout this subsection, assume $l_i \in \mathbb{Z}_{\geq 0}$ for $i=0,1,2,3$.
Assume $\tilde{\beta}_i\in \mathbb{Z}$ $(i=0,1,2,3)$ and
$-\sum_{i=0}^3 \tilde{\beta}_i /2 \in \mathbb{Z}_{\geq 0}$.
Let $V_{\tilde{\beta}_0, \tilde{\beta}_1, \tilde{\beta}_2, \tilde{\beta}_3}$ be a
vector space spanned by
$\big\{ (\frac{\sigma_1(x)}{\sigma(x)})^{\tilde{\beta }_1}
(\frac{\sigma_2(x)}{\sigma(x)})^{\tilde{\beta }_2}
(\frac{\sigma_3(x)}{\sigma(x)})^{\tilde{\beta }_3}\wp(x)^n\big\}
_{n=0, \dots ,-\sum_{i=0}^3 \tilde{\beta}_i /2 }$
and we set $\tilde{V}_{\tilde{\beta}_0, \tilde{\beta}_1, \tilde{\beta}_2,
\tilde{\beta}_3}= \{ \frac{f(x)}{\Phi(x)} | f(x) \in V_{\tilde{\beta}_0,
 \tilde{\beta}_1, \tilde{\beta}_2, \tilde{\beta}_3 }\} $, where $\Phi(x)$
was defined in section \ref{sect:trig}. Let  $\tilde{\alpha}_i \in \{ -l_i , l_i+1\}$
 $(i=0,1,2,3)$ and
\begin{gather*}
 U_{\tilde{\alpha}_0, \tilde{\alpha}_1, \tilde{\alpha}_2, \tilde{\alpha}_3}=
\begin{cases}
V_{\tilde{\alpha}_0, \tilde{\alpha}_1, \tilde{\alpha}_2, \tilde{\alpha}_3},
& \sum_{i=0}^3 \tilde{\alpha}_i/2 \in \mathbb{Z}_{\leq 0} ;\\
V_{1-\tilde{\alpha}_0, 1-\tilde{\alpha}_1, 1-\tilde{\alpha}_2, 1-\tilde{\alpha}_3} ,
& \sum_{i=0}^3 \tilde{\alpha}_i /2\in \mathbb{Z}_{\geq 2} ;\\
\{ 0 \} ,& \mbox{otherwise},
\end{cases} \\
 \tilde{U}_{\tilde{\alpha}_0, \tilde{\alpha}_1, \tilde{\alpha}_2,
\tilde{\alpha}_3}=
\begin{cases}
\tilde{V}_{\tilde{\alpha}_0, \tilde{\alpha}_1, \tilde{\alpha}_2, \tilde{\alpha}_3} ,
&  \sum_{i=0}^3 \tilde{\alpha}_i /2\in \mathbb{Z}_{\leq 0}; \\
\tilde{V}_{1-\tilde{\alpha}_0, 1-\tilde{\alpha}_1, 1-\tilde{\alpha}_2,
1-\tilde{\alpha}_3} ,& \sum_{i=0}^3 \tilde{\alpha}_i /2\in \mathbb{Z}_{\geq 2};\\
\{ 0 \} ,& \mbox{otherwise}.
\end{cases}
\end{gather*}
If $l_0 +l_1 +l_2 +l_3$ is even, then the Hamiltonian $H$ (see (\ref{Ino}))
preserves the spaces \\
$U_{-l_0,-l_1,-l_2,-l_3}$, $U_{-l_0 ,-l_1,l_2+1 ,l_3+1}$,
$U_{-l_0,l_1+1,-l_2,l_3+1}$, $U_{-l_0 ,l_1+1,l_2+1,-l_3}$. The \ 
gauge-transformed Hamiltonian $H(p)$ (see (\ref{eqn:gtrans})) preserves the spaces
$\tilde{U}_{-l_0,-l_1,-l_2,-l_3}$, $\tilde{U}_{-l_0 ,-l_1,l_2+1 ,l_3+1}$,
 $\tilde{U}_{-l_0,l_1+1,-l_2,l_3+1}$, 
$\tilde{U}_{-l_0 ,l_1+1,l_2+1,-l_3}$.
If  $l_0 +l_1 +l_2 +l_3$ is odd, then the Hamiltonian $H$ preserves the spaces
$U_{-l_0,-l_1,-l_2,l_3+1}$, $U_{-l_0 ,-l_1,l_2+1,-l_3}$, $U_{-l_0 ,l_1+1,-l_2,-l_3}$,
 $U_{l_0+1,-l_1,-l_2,-l_3}$, and the gauge-transformed Hamiltonian $H(p)$ preserves
the spaces $\tilde{U}_{-l_0,-l_1,-l_2,l_3+1}$, $\tilde{U}_{-l_0 ,-l_1,l_2+1,-l_3}$,
$\tilde{U}_{-l_0 ,l_1+1,-l_2,-l_3}$, and \break
 $\tilde{U}_{l_0+1,-l_1,-l_2,-l_3}$.

Now we present results on the inclusion of
$ \tilde{U}_{\tilde{\alpha}_0, \tilde{\alpha}_1, \tilde{\alpha}_2, \tilde{\alpha}_3}$
in a Hilbert space. We consider eight cases.

\subsubsection{The case $l_0>0 ,\: l_1>0$ and $l_0+l_1+l_2+l_3$ is even}
If $l_2+l_3-l_0-l_1\geq 2$, then $\tilde{V}_{l_0+1, l_1+1,-l_2, -l_3} \subset
\mathbf{H}$ and $\dim \tilde{V}_{l_0+1, l_1+1,-l_2, -l_3} 
= (l_2+l_3-l_0-l_1)/2$.

\subsubsection{The case $l_0>0 ,\: l_1>0$ and $l_0+l_1+l_2+l_3$ is odd}
If $l_2-l_0-l_1-l_3\geq 3$, then $\tilde{V}_{l_0+1, l_1+1,-l_2, l_3+1} \subset
 \mathbf{H}$ and $\dim \tilde{V}_{l_0+1, l_1+1,-l_2, l_3+1}
= (l_2-l_0-l_1-l_3-1)/2$. If $l_3-l_0-l_1-l_2\geq 3$, then
$\tilde{V}_{l_0+1, l_1+1,l_2+1, -l_3} \subset \mathbf{H}$ and
$\dim \tilde{V}_{l_0+1, l_1+1,l_2+1, -l_3} = (l_3-l_0-l_1-l_2-1)/2$.

\subsubsection{The case $l_0>0 , \: l_1=0$ and $l_0+l_1+l_2+l_3$ is even}
If $l_2+l_3-l_0\geq 2$, then the space $\tilde{V}_{l_0+1, 1,-l_2, -l_3}$ is a
subspace of the space $\mathbf{H}_- (\subset \mathbf{H}) $ and
$\dim \tilde{V}_{l_0+1, 1,-l_2, -l_3} = (l_2+l_3-l_0)/2$.
If $l_2-l_0-l_3\geq 2$, then $\tilde{V}_{l_0+1, 0,-l_2, l_3+1} \subset \mathbf{H}_+$
and $\dim \tilde{V}_{l_0+1, 0,-l_2, l_3+1} = (l_2-l_0-l_3)/2$.
If $l_3-l_0-l_2\geq 2$, then $\tilde{V}_{l_0+1, 0,l_2+1, -l_3} \subset \mathbf{H}_+$
and $\dim \tilde{V}_{l_0+1, 0,l_2+1, -l_3} = (l_3-l_0-l_2)/2$.

\subsubsection{The case $l_0>0 ,\: l_1=0$ and $l_0+l_1+l_2+l_3$ is odd}
If $l_2+l_3-l_0\geq 1$, then $\tilde{V}_{l_0+1, 0,-l_2, -l_3} \subset \mathbf{H}_+$
and $\dim \tilde{V}_{l_0+1, 0,-l_2, -l_3} = (l_2+l_3-l_0+1)/2$.
If $l_2-l_0-l_3\geq 3$, then $\tilde{V}_{l_0+1, 1,-l_2, l_3+1} \subset \mathbf{H}_-$
and $\dim \tilde{V}_{l_0+1, 1,-l_2, l_3+1} =(l_2-l_0-l_3-1)/2$.
If $l_3-l_0-l_2\geq 3$, then $\tilde{V}_{l_0+1, 1,l_2+1, -l_3} \subset \mathbf{H}_-$
and $\dim \tilde{V}_{l_0+1, 1,l_2+1, -l_3} = (l_3-l_0-l_2-1)/2$.

\subsubsection{The case $l_0=0 ,\: l_1>0$ and $l_0+l_1+l_2+l_3$ is even}
If $l_2+l_3-l_1\geq 2$, then $\tilde{V}_{ 1,l_1+1, -l_2, -l_3} \subset \mathbf{H}_-$
and $\dim \tilde{V}_{1, l_1+1, -l_2, -l_3} = (l_2+l_3-l_1)/2$.
If $l_2-l_1-l_3\geq 2$, then $\tilde{V}_{0, l_1+1, -l_2, l_3+1} \subset \mathbf{H}_+$
and $\dim \tilde{V}_{0, l_1+1, -l_2, l_3+1} = (l_2-l_1-l_3)/2$.
If $l_3-l_1-l_2\geq 2$, then $\tilde{V}_{0, l_1+1, l_2+1, -l_3} \subset \mathbf{H}_+$
and $\dim \tilde{V}_{0, l_1+1, l_2+1, -l_3} = (l_3-l_1-l_2)/2$.

\subsubsection{The case $l_0=0 ,\: l_1>0$ and $l_0+l_1+l_2+l_3$ is odd}
If $l_2+l_3-l_1\geq 1$, then $\tilde{V}_{0, l_1+1, -l_2, -l_3} \subset \mathbf{H}_+$
and $\dim \tilde{V}_{0, l_1+1,-l_2, -l_3} = (l_2+l_3-l_1+1)/2$.
If $l_2-l_1-l_3\geq 3$, then $\tilde{V}_{1, l_1+1,-l_2, l_3+1} \subset \mathbf{H}_-$
and $\dim \tilde{V}_{1, l_1+1,-l_2, l_3+1} = (l_2-l_1-l_3-1)/2$.
If $l_3-l_1-l_2\geq 3$, then $\tilde{V}_{1, l_1+1, l_2+1, -l_3} \subset \mathbf{H}_-$
and $\dim \tilde{V}_{1, l_1+1,l_2+1, -l_3} = (l_3-l_1-l_2-1)/2$.

\subsubsection{The case $l_0=l_1=0$ and $l_0+l_1+l_2+l_3$ is even}
In this case, the spaces $\tilde{V}_{0, 0, -l_2, -l_3}, \quad \tilde{V}_{1, 1, -l_2, -l_3},$
\begin{equation*}
\left\{ 
\begin{array}{ll}
\tilde{V}_{1, 0, -l_2, l_3+1},  \quad
\tilde{V}_{0, 1, -l_2, l_3+1} & (l_2>l_3) \\ 
\tilde{V}_{1, 0, l_2+1, -l_3},  \quad 
\tilde{V}_{0,1, l_2+1, -l_3} & (l_2<l_3) \\ 
\{ 0\} & (l_2=l_3)
\end{array}
\right.
\end{equation*}
are subspaces of the Hilbert space $\mathbf{H}$.

\subsubsection{The case $l_0=l_1=0$ and $l_0+l_1+l_2+l_3$ is odd}
In this case, the spaces $\tilde{V}_{1, 0, -l_2, -l_3}, \quad \tilde{V}_{0, 1, -l_2, -l_3},$
\begin{equation*}
\left\{ 
\begin{array}{ll}
\tilde{V}_{1, 1, -l_2, l_3+1}, \quad
\tilde{V}_{0, 0, -l_2, l_3+1} &  (l_2>l_3) \\ 
\tilde{V}_{1, 1, l_2+1, -l_3}, \quad 
\tilde{V}_{0,0, l_2+1, -l_3} & (l_2<l_3)   \\
\{ 0\} & (l_2=l_3)
\end{array}
\right.
\end{equation*}
are subspaces of the Hilbert space $\mathbf{H}$.

\section{Examples} \label{sec:example}

In this section, we show examples that illustrate results of this paper.

\subsection{The case $l_0=1, l_1=2, l_2=0, l_3=8$}
In this case, the Hamiltonian is
\begin{equation*}
H= -\frac{d^2}{dx^2} + 2\wp (x)+6\wp (x+1/2)+72\wp (x+\tau/2).
\end{equation*}
Set $\Phi (x)=(\sin \pi x )^2(\cos \pi x)^3$ and consider the gauge-transformed
operator $H(p)= \Phi (x)^{-1} \circ H \circ \Phi (x)$, where
$p=\exp (\pi\sqrt{-1}\tau )$. By the trigonometric limit $p\to 0$, we obtain
$H(0) =\mathcal{H}_T-\frac{80\pi^2}{3}$, where
\begin{equation}
\mathcal{H}_T= -\frac{d^2}{dx^2}-2\pi\big( \frac{2\cos \pi x}
{\sin \pi x}- \frac{3\sin \pi x}{\cos \pi x} \big) \frac{d}{dx}+25\pi^2.
\label{JacobiHam1208}
\end{equation}

The Hilbert space $\mathbf{H}$ and the inner product are defined by (\ref{innerprod},
\ref{Hilb1}). Then the eigenfunction of the gauge-transformed trigonometric
Hamiltonian $H(0)$ with the eigenvalue $(2m+5)^2\pi^2 -\frac{80 \pi^2}{3}$ is the
Jacobi polynomial $\psi_m(x) = \psi^{(1,2)}_m(x) $ (see (\ref{Jacobipol})).
The functions $\langle \psi_m(x)  \rangle _{m \in \mathbb{Z}_{\geq 0}}$ form an
orthogonal system on the Hilbert space $\mathbf{H}$.
Let $\tilde{H}(p)$ be the selfadjoint extension of $H(p)$.
 From Theorem \ref{mainthmKato} and Corollary \ref{cor:ineq}, it is shown that all
eigenvalues of  $\tilde{H}(p)$ $(-1<p<1)$ on the Hilbert space $\mathbf{H}$ can be
represented as $\tilde{E}_m(p)$ $(m\in \mathbb{Z}_{\geq 0})$, which is real-holomorphic
in $p \in (-1,1)$ and $\tilde{E}_m(0) =(2m+5)^2\pi^2 -\frac{80 \pi^2}{3}$.
Moreover, $\tilde{E}_m(p) < \tilde{E}_{m'}(p)$ for $m<m'$ and $-1<p<1$.

Define spaces of doubly periodic functions as follows:
\begin{gather*}
 V_{2,3,1,-8} =\Big\langle \big(\frac{\sigma_1(x)}{\sigma(x)}\big)^{3
 }\big(\frac{\sigma_2(x)}{\sigma(x)}\big)
 \big(\frac{\sigma_3(x)}{\sigma(x)}\big)^{-8}\wp(x)^n\Big\rangle _{n=0, 1}, \\
 V_{-1,-2,1,-8} =\Big\langle \big(\frac{\sigma_1(x)}{\sigma(x)}\big)^{-2}
 \big(\frac{\sigma_2(x)}{\sigma(x)}\big)\big(\frac{\sigma_3(x)}{\sigma(x)}\big)^{-8}
 \wp(x)^n\Big\rangle _{n=0, \dots, 5}, \nonumber \\
 V_{2,-2,0,-8} =\Big\langle \big(\frac{\sigma_1(x)}{\sigma(x)}\big)^{-2}
 \big(\frac{\sigma_3(x)}{\sigma(x)}\big)^{-8}\wp(x)^n\Big\rangle _{n=0, \dots ,4},
  \nonumber \\
 V_{-1,3,0,-8} =\Big\langle \big(\frac{\sigma_1(x)}{\sigma(x)}\big)^{3}
 \big(\frac{\sigma_3(x)}{\sigma(x)}\big)^{-8}\wp(x)^n\Big\rangle _{n=0, \dots ,3},
 \nonumber
\end{gather*}
and set $\tilde{V}_{\tilde{\beta}_0, \tilde{\beta}_1,
\tilde{\beta}_2, \tilde{\beta}_3}= \{ \frac{f(x)}{\Phi(x)} :
f(x) \in V_{\tilde{\beta}_0, \tilde{\beta}_1, \tilde{\beta}_2, \tilde{\beta}_3 }\}$.
Then the Hamiltonian $H$ preserves  the spaces
$ V_{2,3,1,-8}$, $V_{-1,-2,1,-8}$, $V_{2,-2,0,-8}$,  $V_{-1,3,0,-8} $ and the
gauge-transformed Hamiltonian $H(p)$ preserves the spaces 
$ \tilde{V}_{2,3,1,-8}$,
$\tilde{V}_{-1,-2,1,-8}$, $\tilde{V}_{2,-2,0,-8}$, $\tilde{V}_{-1,3,0,-8} $.
Among the spaces $ \tilde{V}_{2,3,1,-8}$, $\tilde{V}_{-1,-2,1,-8}$,
$\tilde{V}_{2,-2,0,-8} $, $\tilde{V}_{-1,3,0,-8} $, only the space $\tilde{V}_{2,3,1,-8}$ is a
subspace of the Hilbert space $\mathbf{H}$. The eigenvalues of the gauge-transformed
Hamiltonian $H(p)$ on the space $\tilde{V}_{2,3,1,-8}$ are written as
$11e_1-9e_2 \pm 2\sqrt{106e_1^2+73e_1e_2+46e_2^2}$, where $e_1=\wp(1/2)$ and
$e_2=\wp((1+\tau)/2)$. From Theorem \ref{thm:lowest}, we obtain
$\tilde{E}_0(p)= 11e_1-9e_2 - 2\sqrt{106e_1^2+73e_1e_2+46e_2^2}$ and
$\tilde{E}_1(p)= 11e_1-9e_2 + 2\sqrt{106e_1^2+73e_1e_2+46e_2^2}$.
Hence the smallest eigenvalue and the second smallest one on the Hilbert space
$\mathbf{H}$ are obtained algebraically.

\subsection{The case $l_0=1, l_1=2, l_2=1, l_3=0$}
In this case, the Hamiltonian is
\begin{equation*}
H= -\frac{d^2}{dx^2} + 2\wp (x)+6\wp (x+1/2)+2\wp (x+(\tau+1)/2).
\end{equation*}
Set $\Phi (x)=(\sin \pi x )^2(\cos \pi x)^3$ and consider the gauge-transformation
$H(p)= \Phi (x)^{-1} \circ H \circ \Phi (x)$.
The Hilbert space $\mathbf{H}$ is defined similarly to the case
$l_0=1, l_1=2, l_2=0, l_3=8$.
Define
\begin{gather*}
 V_{-1,-2,-1,0} =\Big\langle \big(\frac{\sigma_1(x)}{\sigma(x)}\big)^{-2}
 \big(\frac{\sigma_2(x)}{\sigma(x)}\big)^{-1}\wp(x)^n\Big\rangle _{n=0, 1,2}, \\
 V_{-1,-2,2,1} =\Big\langle \big(\frac{\sigma_1(x)}{\sigma(x)}\big)^{-2}
 \big(\frac{\sigma_2(x)}{\sigma(x)}\big)^2 \big(\frac{\sigma_3(x)}{\sigma(x)}\big)
 \Big\rangle , \nonumber \\
V_{2,-2,-1,1} =\Big\langle \big(\frac{\sigma_1(x)}{\sigma(x)}\big)^{-2}
\big(\frac{\sigma_2(x)}{\sigma(x)}\big)^{-1}\big(\frac{\sigma_3(x)}{\sigma(x)}\big)
\Big\rangle , \nonumber
\end{gather*}
and set $\tilde{V}_{\tilde{\beta}_0, \tilde{\beta}_1, \tilde{\beta}_2,
\tilde{\beta}_3}= \big\{ \frac{f(x)}{\Phi(x)} | f(x) \in V_{\tilde{\beta}_0,
\tilde{\beta}_1, \tilde{\beta}_2, \tilde{\beta}_3 }\big\}$.
Then the Hamiltonian $H$ preserves the spaces $ V_{-1,-2,-1,0}$, $V_{-1,-2,2,1}$,
$V_{2,-2,-1,1} $ and the gauge-transformed \break
Hamiltonian $H(p)$ preserves the spaces
$ \tilde{V}_{-1,-2,-1,0}$, $\tilde{V}_{-1,-2,2,1}$, $\tilde{V}_{2,-2,-1,1} $.
But none of spaces $ \tilde{V}_{-1,-2,-1,0}$, $\tilde{V}_{-1,-2,2,1}$,
$\tilde{V}_{2,-2,-1,1} $ is included in the Hilbert space $\mathbf{H}$.

\subsection{The case $l_0=1, l_1=0, l_2=4, l_3=1$}
In this case, the Hamiltonian is
\begin{equation*}
H= -\frac{d^2}{dx^2} + 2\wp (x)+20\wp (x+(\tau+1)/2)+2\wp (x+\tau/2).
\end{equation*}
Set $\Phi (x)=\sin \pi x $ and consider the gauge-transformation $H(p)= \Phi (x)^{-1}\circ H \circ \Phi (x)$, where $p=\exp(\pi\sqrt{-1}\tau )$.
By the trigonometric limit $p\to 0$, we obtain $H(0) =\mathcal{H}_T-8\pi^2$, where
\begin{equation}
\mathcal{H}_T= -\frac{d^2}{dx^2}-2\pi\left( \frac{2\cos \pi x}{\sin \pi x}- \frac{3\sin \pi x}{\cos \pi x} \right) \frac{d}{dx}+4\pi^2.
\label{JacobiHam1041}
\end{equation}

The Hilbert space $\mathbf{H}$, its subspaces $\mathbf{H}_+$, $\mathbf{H}_-$, and the
inner product are defined by (\ref{innerprod}, \ref{Hilbge}). Then the eigenfunction
of the gauge-transformed trigonometric Hamiltonian $H(0)$ with the eigenvalue
$(m+2)^2\pi^2 -8 \pi^2$ is written by use of the Gegenbauer polynomial
$\psi^G_m(x)$ (see (\ref{Gegenpol})).
The functions $\langle \psi^G_m(x)  \rangle _{m \in \mathbb{Z}_{\geq 0}}$
(resp.  $\langle \psi^G_m(x)  \rangle _{m \in 2\mathbb{Z}_{\geq 0}}$,
$\langle \psi^G_m(x)  \rangle _{m \in 2\mathbb{Z}_{\geq 0}+1}$) form an orthogonal
system on the Hilbert space $\mathbf{H}$ (resp. $\mathbf{H}_+$, $\mathbf{H}_-$).
Let $\tilde{H}(p)$ be the selfadjoint extension of $H(p)$.
From Theorem \ref{mainthmKato} and Corollary \ref{cor:ineq}, it is shown that all
eigenvalues of  $\tilde{H}(p)$ $(-1<p<1)$ on the Hilbert space $\mathbf{H}$ can be
represented as $\tilde{E}_m(p)$, which is real-holomorphic in $p \in (-1,1)$ and
$\tilde{E}_m(0) =(m+2)^2\pi^2 -8 \pi^2$. Define
\begin{gather*}
 V_{-1,0,-4,-1} =\Big\langle \big(\frac{\sigma_2(x)}{\sigma(x)}\big)^{-4}
 \big(\frac{\sigma_3(x)}{\sigma(x)}\big)^{-1}\wp(x)^n\Big\rangle _{n=0, 1,2,3}, \\
 V_{2,1,-4,-1} =\Big\langle \big(\frac{\sigma_1(x)}{\sigma(x)}\big)
 \big(\frac{\sigma_2(x)}{\sigma(x)}\big)^{-4} \big(\frac{\sigma_3(x)}
 {\sigma(x)}\big)^{-1}\wp(x)^n\Big\rangle _{n=0, 1}, \nonumber \\
 V_{2,0,-4,2} =\Big\langle \big(\frac{\sigma_2(x)}{\sigma(x)}\big)^{-4}
 \big(\frac{\sigma_3(x)}{\sigma(x)}\big)^{2}\big\rangle, \nonumber \\
 V_{-1,1,-4,2} =\Big\langle \big(\frac{\sigma_1(x)}{\sigma(x)}\big)
 \big(\frac{\sigma_2(x)}{\sigma(x)}\big)^{-4}\big(\frac{\sigma_3(x)}
 {\sigma(x)}\big)^{2}\wp(x)^n\Big\rangle _{n=0, 1}, \nonumber
\end{gather*}
and set $\tilde{V}_{\tilde{\beta}_0, \tilde{\beta}_1, \tilde{\beta}_2,
 \tilde{\beta}_3}= \{ \frac{f(x)}{\Phi(x)} : f(x) \in V_{\tilde{\beta}_0,
 \tilde{\beta}_1, \tilde{\beta}_2, \tilde{\beta}_3 }\}$.
Then the Hamiltonian $H$ preserves spaces $ V_{-1,0,-4,-1}$, $V_{2,1,-4,-1}$,
$V_{2,0,-4,2} $, $V_{-1,1,-4,2} $, the gauge-transformed Hamiltonian $H(p)$
preserves spaces $ \tilde{V}_{-1,0,-4,-1}$, $\tilde{V}_{2,1,-4,-1}$,
$\tilde{V}_{2,0,-4,2} $, $\tilde{V}_{-1,1,-4,2} $, and we have
$\tilde{V}_{2,0,-4,2} \subset \mathbf{H}_+$, and $\tilde{V}_{2,1,-4,-1}\subset
\mathbf{H}_-$. 

The eigenvalue of the gauge-transformed Hamiltonian $H(p)$ on the
space $\tilde{V}_{2,0,-4,2}$ is $-12e_1-12e_2$ and the eigenvalues of $H(p)$ on
the space $\tilde{V}_{2,1,-4,-1}$ are written as
$12e_1-3e_2 \pm 2\sqrt{39e_1^2+3e_1e_2-6e_2^2}$.
 From the results in section \ref{sec:L2l0l12}, the smallest eigenvalue on
$\mathbf{H}_+$ is $-12e_1-12e_2$ and the smallest two eigenvalues on $\mathbf{H}_-$
are $12e_1-3e_2 \pm 2\sqrt{39e_1^2+3e_1e_2-6e_2^2}$. Hence
$\tilde{E}_0(p)= -12e_1-12e_2$, $\tilde{E}_1(p)=12e_1-3e_2 -2\sqrt{39e_1^2
+3e_1e_2-6e_2^2}$, and $\tilde{E}_3(p)=12e_1-3e_2 + 2\sqrt{39e_1^2+3e_1e_2-6e_2^2}$.
In other words, the smallest eigenvalue, the second smallest one, and the fourth
smallest one on the Hilbert space $\mathbf{H}$ are obtained algebraically.

\section{Concluding remarks} \label{sec:com}

In the paper \cite{GGR}, Gomez-Ullate, Gonzalez-Lopez, and Rodriguez found
square-integrable finite dimensional invariant spaces for the $BC_N$ Inozemtsev
model with some coupling constants, and they observed numerically for an example of
the $BC_2$ Inozemtsev model that the set of eigenvalues on a finite dimensional
 invariant space would coincide with the set of small eigenvalues of the
Hamiltonian on $L^2$ space from the bottom.
In the present paper, we justified this phenomena concretely for the $BC_1$
Inozemtsev model in Theorem \ref{thm:lowest} and sections \ref{sec:L2l0l12},
\ref{sec:L2l0l13}, and we illustrated plainly with examples in section
\ref{sec:example}. It would be possible to obtain similar results for the
$BC_N$ Inozemtsev model by combining the Kato-Rellich theory and the trigonometric
limit.

We comment on relationship to the known facts for the Heun equation.
 In books \cite{Ron,SL}, notions ``Heun polynomial'' and ``Heun function'' are
introduced. Roughly speaking, the Heun polynomial corresponds to a doubly periodic
eigenfunction of the $BC_1$ Inozemtsev model, and the Heun function corresponds to
a square-integrable eigenfunction or a doubly periodic eigenfunction of the
$BC_1$ Inozemtsev model.
We hope that both the Heun equation and the Inozemtsev model are elusidated
substantially in near future.

\section{Appendix} \label{sec:app}

\subsection{Proof of Proposition \ref{prop:pmvm}}
In this subsection we prove Proposition \ref{prop:pmvm}.
\begin{lemma} \label{anallem2}
Let $\psi_m (x)$ be the normalized eigenbasis of the gauge-transformed 
trigonometric Hamiltonian $\mathcal{H}_T$ (\ref{JacobiHam}) defined in 
(\ref{Jacobipol}), and write $\sum_{k=1}^{\infty} V_k(x) p^k =H(p)-H(0)$. Then
\begin{equation}
\psi_m (x) V_k(x) = \sum_{|{m'}-m|\leq k }\bar{c}_{{m'}} \psi_{m'} (x),
\label{psinpieri}
\end{equation}
for some constants $\bar{c}_{{m'}} $ $(|{m'}-m|\leq k )$.
\end{lemma}
\begin{proof}
 From the Pieri-type formula for the Jacobi polynomials, it follows that
\begin{equation}
\psi_m (x) \cos 2\pi x = c_{-}\psi_{m-1} (x) +  c_{0}\psi_{m} (x) +c_{+}\psi_{m+1} (x) ,
\label{psipieri}
\end{equation}
for some constants $c_{-}, c_{0}$ and $c_{+}$.
On the other hand, the function $V_k(x)$ is a polynomial in the variable 
$\cos 2\pi x$ with degree $k$. By applying (\ref{psipieri}) repeatedly, 
we obtain (\ref{psinpieri}).
\end{proof}

\begin{prop} \label{analprop1}
Let $|p|<1$ and  $(\sum_{k=1}^{\infty}V_k(x) p^{k})\psi _m(x)
=\sum_{{m'}\in \mathbb{Z} _{\geq 0}}\tilde{t}_{m.{m'}}\psi_ {m'}(x)$. 
For each $C$ such that  $C>1$ and $C|p|<1$, there exists 
$C''\in \mathbb{R}_{>0}$ such that 
$|\tilde{t}_{m,{m'}}| \leq C''(C|p|)^{\frac{|{m'} -m|+1}{2}}$.
\end{prop}

\begin{proof}
Since the normalized Jacobi polynomials form a complete orthonormal system,
it follows that $\tilde{t}_{m,{m'}}= \langle \psi_ {m'}(x) ,(\sum_{k=1}^{\infty}V_k(x)p^{k})\psi _m(x)\rangle _{\Phi }$.

If $k<|m'-m|$, then $ \langle V_k(x)p^{k}\psi _m(x), \psi_ {m'}(x) \rangle _{\Phi }=0$ by Lemma \ref{anallem2} and orthogonality.
Therefore,
\begin{align*}
 |\tilde{t}_{m,{m'}}|
& = \Big| \Big\langle \psi_ {m'}(x), \Big(\sum_{k=1}^{\infty}V_k(x)p^{k}\Big)
\psi _m(x) \Big\rangle _{\Phi } \Big|\\
& = \Big| \Big\langle  \psi_{m'}(x),\Big(\sum_{k=|m'-m|}^{\infty}V_k(x)p^{k}\Big)
\psi _m(x) \Big\rangle _{\Phi } \Big|  \\
& = \Big|\int _0^1 \sum_{k=|m'-m|}^{\infty}V_k(x)p^{k}\psi _m(x) 
\overline{\psi_ {m'}(x)} \Phi(x) ^2 dx \Big|  \\
& \leq \sup_{x\in [0,1]} \Big| \sum_{k=|m'-m|}^{\infty}V_k(x)p^{k} \Big| 
\int_0^1 |\psi _m(x) \psi_ {m'}(x) \Phi(x)^2 | dx  \\
& \leq \sum_{k \geq |m'-m|}V_k|p|^{k}. 
\end{align*}
where $V_k$ is defined by (\ref{majser}).
For the case $m'=m$, we obtain 
$|\tilde{t}_{m,m}| \leq \sum_{k \geq 1}V_k|p|^{k}$.

Since the convergence radius of the series $\sum_n V_n p^n$ is equal to $1$, 
for each $C$ such that $C>1$ and $C|p|<1$ there exists $C''\in \mathbb{R} _{>0}$ 
such that $|\tilde{t}_{m,{m'}}| \leq C''(C|p|)^{|{m'} -m|}$ $(m\neq {m'})$ and 
$|\tilde{t}_{m,m}| \leq C''(C|p|) \leq C''(C|p|)^{1/2}$ . From an inequality 
$(C|p|)^{|{m'} -m|} \leq (C|p|)^{\frac{|{m'} -m|+1}{2}}$ for 
$|{m'}-m| \in \mathbb{Z}_{\geq 1}$, we obtain the proposition.
\end{proof}

\begin{prop} \label{analprop11}
Let $D$ be a positive number. Suppose 
$\mathop{\rm dist}(\zeta, \sigma (\tilde{H}(0)))\geq D$. Write 
$(\tilde{H}(p)-\zeta )^{-1}  \psi _m(x)=\sum_{{m'}}t_{m,{m'}}\psi_ {m'}(x)$, 
where $(\tilde{H}(p)-\zeta )^{-1}$ is defined in (\ref{Neus}).
For each $m \in \mathbb{Z}_{\geq 0}$ and $C\in \mathbb{R} _{>1}$, there exists 
$C'\in \mathbb{R} _{>0}$ and $p_0\in \mathbb{R} _{>0}$ which do not depend on 
$\zeta $ (but depend on $D$) such that $t_{m,{m'}}$ satisfy
\begin{equation}
|t_{m,{m'}}| \leq C'(C|p|)^{|{m'} -m|/2},
\label{proptlm}
\end{equation}
 for all $p$, $m'$ such that $|p|<p_0$ and ${m'} \in \mathbb{Z}_{\geq 0}$.
\end{prop}
\begin{proof}
Recall that the operator $(\tilde{H}(p)-\zeta )^{-1}$ is defined by Neumann 
series (\ref{Neus}).
Fix a number $D(\in \mathbb{R}_{>0})$ and set 
$X:=(\zeta -\tilde{H}(0))^{-1}(\sum_{k=1}^{\infty}V_k(x)p^{k})$.
 From expansion (\ref{Neus}), there exists $p_1\in \mathbb{R}_{>0}$ such that 
 an inequality $\| X \| <1/2$ holds for all $p$ and $\zeta$ such that $|p|<p_1$ 
and $\mathop{\rm dist}(\zeta, \sigma (\tilde{H}(0)))>D$. Then we have 
$(\sum_{i=0}^{\infty}X^i)(\tilde{H}(0)-\zeta )^{-1}= (\tilde{H}(p)-\zeta )^{-1}$.
If we write $\sum_{{m'}} c'_{{m'}}\psi_ {m'}(x)=(\tilde{H}(0)-\zeta )^{-1} 
\sum_{{m'}}c_{{m'}}\psi_ {m'}(x)$, then $|c'_{{m'}}| \leq D^{-1}|c_{{m'}}|$ 
for each ${m'}$.
Write $X\psi _m(x)=\sum_{{m'}}\check{t}_{m,{m'}}\psi_ {m'}(x)$.
By combining with Proposition \ref{analprop1}, we obtain that for each $C$ 
such that $C>1$ and $\frac{C+1}{2} p_1<1$, there exists $C''\in \mathbb{R} _{>0}$
which does not depend on $\zeta $ (but depend on $D$) such that 
$|\check{t}_{m,{m'}}| \leq C''(\frac{C+1}{2}|p|)^{(|{m'} -m|+1)/2}$ for 
$|p|<p_1$.

To obtain Proposition \ref{analprop11}, we use the method of majorants.
For this purpose, we introduce symbols $\mathbf{e}_{m }$ $(m \in \mathbb{Z} )$ 
to avoid inaccuracies and apply the method of majorants for formal series 
$\sum_{{m'}\in \mathbb{Z}} c_{{m'}}\mathbf{e}_{{m'}}$.
For formal series, we define the partial ordering $\Tilde{\leq }$ by the 
following rule:
\begin{equation*}
\sum_{{m'}\in \mathbb{Z}} c^{(1)}_{{m'}}\mathbf{e}_{{m'}} \Tilde{\leq } 
\sum_{{m'}\in \mathbb{Z}} c^{(2)}_{{m'}}\mathbf{e}_{{m'}}\;
\Leftrightarrow\; \forall {m'}, \;  |c^{(1)}_{{m'}}|\leq |c^{(2)}_{{m'}}| .
\end{equation*}
We will later consider the case that each coefficient $c_{{m'}}^{(i)}$ 
$(i=1,2, {m'} \in \mathbb{Z})$ is expressed as an infinite sum.
If absolute convergence of $c_{{m'}}^{(2)}$ is shown for each ${m'}$, 
then absolute convergence of $c_{{m'}}^{(1)}$ for each ${m'}$ is shown by 
the majorant.
Set $X\mathbf{e}_{m}=\sum_{{m'}}\check{t}_{m,{m'}}\mathbf{e}_{{m'}}$, where 
coefficients $\check{t}_{m,{m'}}$ are defined by 
$X\psi _m(x)=\sum_{{m'}}\check{t}_{m,{m'}}\psi_ {m'}(x)$.

Our goal is to show (\ref{proptlm}) for $t_{m,{m'}}$ such that 
$\sum_{{m'}} t_{m,{m'}} \psi_{{m'} }(x) =  \sum_{i=0}^{\infty} X^i 
(\tilde{H}(0)-\zeta )^{-1}\psi _m(x)$.
Since $(\tilde{H}(0)-\zeta )^{-1}\psi _m(x)= (E_{m} -\zeta )^{-1}\psi _m(x)$ 
and $|(E_{m} -\zeta )^{-1}|\leq D^{-1}$, it is enough to show that
there exist $C^{\star} \in \mathbb{R} _{>0}$ and $p_0$ which do not depend on 
$\zeta $ (but depend on $D$) such that $t^{\star}_{m,{m'}}$ are well-defined 
by  $\sum_{k=0}^{\infty} X^k \mathbf{e}_{m}= \sum_{{m'}\in \mathbb{Z}}
 t^{\star}_{m,{m'}}\mathbf{e}_{{m'}}$ and satisfy
\begin{equation}
|t^{\star}_{m,{m'}}| \leq C^{\star}\Big(\frac{C+1}{2}|p|\Big)^{|{m'} -m|/2},
\label{proptlms}
\end{equation}
 for all $p$ and ${m'}$ such that $|p|<p_0$ and ${m'} \in \mathbb{Z}_{\geq 0}$.
Set
\begin{equation*}
Z \mathbf{e}_{m} = \sum_{m \in \mathbb{Z} }z_{m,{m'}} \mathbf{e}_{{m'}}
=\sum_{{m'} \in \mathbb{Z}}C''\Big(\frac{C+1}{2}p\Big)^{\frac{|m -{m'} |+1}{2}} 
\mathbf{e}_{{m'}}.
\end{equation*}
Then we have $X \mathbf{e}_{m} \Tilde{\leq } Z \mathbf{e}_{m}$.
Let $k \in \mathbb{Z}_{\geq 1}$. If coefficients of $Z^k \mathbf{e}_{m}$ 
with respect to the basis $\{ \mathbf{e}_{{m'}} \}$ converge absolutely, 
then $X^k \mathbf{e}_{m}$ is well-defined and
$X^k \mathbf{e}_{m} \Tilde{\leq } Z^k \mathbf{e}_{m}$.

 From the equality $Z^k \mathbf{e}_{m}= \sum_{\nu^{(1)}, 
\dots, \nu^{(k-1)}} z_{m , \nu^{(1)}}z_{\nu^{(1)},\nu^{(2)}}
\dots z_{\nu^{(k-1)},{m'}} \mathbf{e}_{{m'}}$ and the property 
$z_{m,{m'} }=z_{0,{m'}-m}$, it follows that
\begin{align*}
  \sum_{k=1}^{\infty} Z^k \mathbf{e}_{m}
&= \sum_{k=1}^{\infty}\sum_{{m'} \in \mathbb{Z} }\frac{1}{2\pi\sqrt{-1}}
\oint _{|s|=1} \Big(\sum_{\nu \in  \mathbb{Z} }z_{0,\nu }s^{\nu}\Big)^k 
s^{m-{m'}-1}ds \mathbf{e}_{{m'}}\\
& = \sum_{k=1}^{\infty}\sum_{{m'} \in \mathbb{Z} }\frac{1}{2\pi\sqrt{-1}}
\oint_{|s|=1}\Big(\sum_{\nu \in \mathbb{Z}} C'' \Tilde{p}^{(|\nu |+1)}s^{\nu}
\Big)^k s^{m-{m'}-1}ds \mathbf{e}_{{m'}}  \\
& = \sum_{{m'} \in \mathbb{Z}} Z_{m ,{m'}} \mathbf{e}_{{m'}},  
\end{align*}
where $\Tilde{p}=(\frac{C+1}{2}p)^{1/2}$ and
\begin{equation}
 Z_{m,{m'}}= \frac{1}{2\pi\sqrt{-1}} \oint_{|s|=1}\frac{C''(\tilde{p}
-\tilde{p}^3)s^{m-{m'}-1}ds}{(1-\Tilde{p}s)(1-\Tilde{p}s^{-1})
-C''(\tilde{p}-\tilde{p}^3)} .
\label{zlm}
\end{equation}
Note that we used a formula 
$\sum_{n\in \mathbb{Z}}q^{|n|+1}x^n=\frac{q-q^3}{(1-qx)(1-qx^{-1})}$. 
Equality (\ref{zlm}) make sense for $\tilde{p}<p_2$, where $p_2$ is a positive 
number satisfying inequalities $p_2<1$, $\frac{C''|p_2-p_2^3|}{(1-p_2)^2}<1$ 
and $C''p_2<1$.
Therefore each coefficient of $\sum_{k=1}^{\infty} Z^k \mathbf{e}_{m}$ with 
respect ot the basis $\{ \mathbf{e}_{{m'}} \}$ converges absolutely. 
It follows that
\begin{equation}
\sum_{k=0}^{\infty} X^k \mathbf{e}_{m}  \Tilde{\leq }  \mathbf{e}_{m}
+ \sum_{k=1}^{\infty} Z^k \mathbf{e}_{m} \Tilde{\leq } \mathbf{e}_{m}+ 
 \sum_{{m'}\in \mathbb{Z}} Z_{m,{m'}} \mathbf{e}_{{m'}}.
\end{equation}

Let $s(\Tilde{p})$ be the solution of an equation 
$(1-\Tilde{p}s)(1-\Tilde{p}s^{-1})-C''(\tilde{p}-\tilde{p}^3)=0$ on the 
variable $s$ satisfying $|s(\Tilde{p})|<1$. Then $s(\Tilde{p})$ is 
holomorphic in $\Tilde{p}$ near $0$ and admits an expansion 
$s(\Tilde{p})= \Tilde{p}+ c_2\Tilde{p}^2+\dots $, and
\begin{equation}
\frac{1}{(2\pi\sqrt{-1})}\oint_{|s|=1}\frac{C''(\tilde{p}-\tilde{p}^3)
s^{n-1}ds}{(1-\Tilde{p}s)(1-\Tilde{p}s^{-1})-C''(\tilde{p}-\tilde{p}^3)}
= \Tilde{p} f(\Tilde{p})s(\Tilde{p})^{|n|},
\label{integTp}
\end{equation}
where $f(\Tilde{p})$ is a holomorphic function defined near $\Tilde{p}=0$.
Note that relation (\ref{integTp}) is shown by calculating the residue around 
$s=s(\Tilde{p})$ for the case $n\geq 0$ and we need to change a variable 
$s \to s^{-1}$ and calculate the residue around $s=s(\Tilde{p})$ for the case 
$n<0$.
By combining (\ref{zlm}--\ref{integTp}) there exists a positive number such that
\begin{equation}
\sum_{k=0}^{\infty}X^k \mathbf{e}_{m} \Tilde{\leq }\mathbf{e}_{m} 
+ \sum_{{m'}\in \mathbb{Z}} \Tilde{p} f(\Tilde{p})s(\Tilde{p})^{|{m'} -m |} 
\mathbf{e}_{{m'}}.
\label{emieq}
\end{equation}
for $|\Tilde{p}|<p_3$.
By combining (\ref{emieq}), a relation 
$\Tilde{p}=\left(\frac{C+1}{2}p\right)^{\frac{1}{2}}$, an inequality 
$\frac{C+1}{2}<C$ and an expansion 
$s(\Tilde{p})= \Tilde{p}+ c_2\Tilde{p}^2+\dots $, we obtain (\ref{proptlms}) 
and the proposition.
\end{proof}

The following proposition is essentially the same as 
Proposition \ref{prop:pmvm}.

\begin{prop} \label{analprop2} 
Let $E_m \in \sigma (\tilde{H}(0))$ and $\Gamma _m$ be a circle which 
contains only one element $E_m$ of the set $\sigma (\tilde{H}(0))$ inside it. 
Let $\psi _m(x)$ be the corresponding normalized eigenfunction.
Set $P_m(p)= -\frac{1}{2\pi \sqrt{-1}}\int_{\Gamma_m}(\tilde{H}(p)
-\zeta )^{-1} d\zeta $ and write 
$P_m(p) \psi _m(x)=\sum_{{m'}}s_{m,{m'}}\psi_ {m'}(x)$.
For each $C\in \mathbb{R} _{>1}$, there exist $C'\in \mathbb{R} _{>0}$ 
and $p_{\ast }\in \mathbb{R} _{>0}$ such that $s_{m,{m'}}$ satisfy
\begin{equation}
|s_{m,{m'}}| \leq C'(C|p|)^{|{m'} -m|/2},
\label{propslm}
\end{equation}
 for all $p$ and $m'$ such that $|p|<p_{\ast }$ and $m' \in \mathbb{Z}_{\geq 0}$.
\end{prop}

\begin{proof}
Since the spectrum $\sigma (\Tilde{H}(0))$ is discrete,
there exists a positive number $D$ such that 
$\mathop{\rm dist}(\zeta, \sigma (\Tilde{H}(0))) \geq D$ for 
$\zeta \in \Gamma_m$. Write 
$(\Tilde{H}(p)-\zeta )^{-1} \psi_m(x)=\sum_{{m'}}t_{m,{m'}}(\zeta )\psi_{m'}(x)$.
 From Proposition \ref{analprop11}, it is shown that for each 
 $C\in \mathbb{R} _{>1}$ there exists $C_{\ast }\in \mathbb{R} _{>0}$ and 
 $p_{\ast }\in \mathbb{R} _{>0}$ which do not depend on $\zeta (\in \Gamma_m) $ 
such that $t_{m,{m'}}$ satisfy  $|t_{m,{m'}}(\zeta )| \leq C_{\ast }
(C|p|)^{\frac{|{m'} -m|}{2}}$ for all $p$, $m'$ such that $|p|<p_{\ast }$ and 
$m' \in \mathbb{Z}_{\geq 0}$.

Let $L$ be a length of the circle $\Gamma_m$ and write 
$-\frac{1}{2\pi \sqrt{-1}} \int_{\Gamma_m}(\Tilde{H}(p)-\zeta )^{-1} 
d\zeta \psi_m(x)=\sum_{{m'}}s_{m,{m'}}\psi_{m'}(x)$.
By integrating the function $\sum_{{m'}}t_{m,{m'}}(\zeta )\psi_{m'}(x)$ over 
the circle $ \Gamma_m$ in the variable $\zeta $, it follows that 
$|s_{m,{m'}}| \leq \frac{L}{2\pi} C_{\ast }(C|p|)^{\frac{|{m'} -m|}{2}}$ 
for all $p$, $m'$ such that $|p|<p_{\ast }$ and $m' \in \mathbb{Z}_{\geq 0}$.
Therefore, Proposition \ref{analprop2} is proved.
\end{proof}

\subsection{Proof of Proposition \ref{prop:psival}}
We prove Proposition \ref{prop:psival}.
\begin{prop}[Proposition \ref{prop:psival}] \label{prop:psival2} 
Let $\psi _m(x)$ be the $(m+1)$st normalized eigenfunction of the trigonometric
gauge-transformed Hamiltonian $\mathcal{H}_T$.
Let $f(x)=\sum_{m=0}^{\infty}c_m \psi _m(x)$ be a function satisfying
$|c_m|<A R^m$ $(\forall m \in \mathbb{Z} _{\geq 0})$ for some 
$A \in \mathbb{R}_{>0}$ and $R \in (0,1)$.
If $r'$ satisfies $0<r'<\frac{1}{2\pi} \log \frac{1}{R}$, then the power series 
$\sum_{m=0}^{\infty}c_m \psi _m(x)$ converges uniformly absolutely inside a 
zone $-r' \leq \Im x\leq r'$, where $\Im x$ is an imaginary part of the complex 
number $x$.
\end{prop}

\begin{proof}
We prove for the case $l_0>0$ and $l_1>0$. For the other cases, they are 
proved similarly.
We introduce a Rodrigues-type formula for the Jacobi polynomials
\begin{equation*}
p_m(w)= \frac{(-1)^m}{m!}w^{-l_0-1/2}(1-w)^{-l_1-1/2}\big( \frac{d}{dw} \big) ^m 
\left( w^{l_0+m+1/2}(1-w)^{l_1+m+1/2} \right).
\end{equation*}
Then $p_m(\sin ^2 \pi x)=d_m \psi _m(x)$, where 
$d_m=\sqrt{\frac{\Gamma (m+l_0+3/2)\Gamma (m+l_1+3/2)}{\pi m! 
(2m+l_0+l_1+2)\Gamma (m+l_0+l_1+2)}}$.
 From the Stirling's formula, we have $(d_m)^{1/m}\to 1$ as $m \to \infty $.
Hence it is sufficient to show that the power series 
$\sum_{m=0}^{\infty}c_m p_m(\sin ^2 \pi x)$ converges uniformly absolutely 
inside the zone $-r' \leq \Im x\leq r'$.

The generating function of the Jacobi polynomials $p_m (w)$ is written as 
\begin{equation}
\sum_{m=0}^{\infty} p_m(w)\xi^m= \frac{1}{S\big( \frac{1+\xi+S}{2} \big) ^{l_0+1/2}
 \big( \frac{1-\xi+S}{2} \big) ^{l_1+1/2}},
\label{genfunct}
\end{equation}
where $S=\sqrt{(1+\xi)^2-4\xi w}$.
Now we set $\sum_{m=0}^{\infty} \tilde{q}_m(y)\xi^m
= \frac{1}{\sqrt{(1-\xi)^2-4y\xi}}$, \\
$\sum_{m=0}^{\infty} \tilde{q}^{(a)}_m(y)\xi^m
= \frac{1}{(1-\xi +\sqrt{(1-\xi)^2-4y\xi})^a}$, and
\begin{equation}
\sum_{m=0}^{\infty} \tilde{p}_m(y)\xi^m =  \frac{1}{(\sqrt{(1-\xi)^2-4y\xi})}
\frac{1}{(1-\xi +\sqrt{(1-\xi)^2-4y\xi})^{l_0+l_1+1}}.
\label{genpny}
\end{equation}
Then it is shown that, if $a>0$, then $\tilde{q}_m(y)$ and 
$\tilde{q}^{(a)}_m(y)$ are polynomials in $y$ of degree $m$ with nonnegative 
coefficients. Hence $\tilde{p}_m(y)$ is also a polynomial in $y$ of degree
$m$ with nonnegative coefficients.
Set $p_m(y)= \sum_{k=0}^m p_m^{(k)} y^k$ and 
$\tilde{p}_m(y)= \sum_{k=0}^m \tilde{p}_m^{(k)} y^k$. From formulas 
(\ref{genfunct}, \ref{genpny}) and the nonnegativity, we obtain  
$|p_m^{(k)}|\leq  \tilde{p}_m^{(k)}$ for all $m$ and $k$.
 From the inequality $|\sin^2\pi x |
 \leq \big| \big(\frac{e^{\pi r'}+e^{-\pi r'}}{2}\big)^2 \big|$ for
$-r'\leq \Im x \leq r'$, it is seen that
\begin{equation}
\begin{aligned}
|p_m(\sin^2\pi x)| &\leq \sum_{k=0}^m \big| p_m^{(k)} ( \sin^2\pi x )^k \big|
\leq \sum_{k=0}^m \tilde{p}_m^{(k)} | \sin^2\pi x |^k \\
&\leq \Big|\tilde{p}_m\Big( \big(\frac{e^{\pi r'}+e^{-\pi r'}}{2}\big)^2\Big)\Big|
\end{aligned}\label{ineq:pn}
\end{equation}
for all $m\in \mathbb{Z} _{\geq 0}$ and $x$ such that $-r'\leq \Im x \leq r'$.

On the other hand,  the series 
$\sum_{m=0}^{\infty} \tilde{p}_m\big ((\frac{e^{\pi r'}+e^{-\pi r'}}{2})^2
\big) \xi ^m$, with respect to the variable $\xi$, has radius of convergence
$e^{-2\pi r'}$; because the singular point of the right hand side of
(\ref{genpny}) which is closest to the origin is located on the circle 
$|\xi |=e^{-2\pi r'}$.
Let $r''$ be a positive number such that 
$r'<r''<\frac{1}{2\pi} \log (1/R)$.
Then there exists $A' \in \mathbb{R} $ such that 
$\tilde{p}_m \big( (\frac{e^{\pi r'}+e^{-\pi r'}}{2})^2\big) < A' e^{2\pi r''m}$.
 Hence we have 
$ |c_m | \big| \tilde{p}_m\big( (\frac{e^{\pi r'}+e^{-\pi r'}}{2})^2\big)\big|
 <A A' \big( R e^{2\pi r''} \big) ^m$.
Since $R e^{2\pi r''} <1$, the series 
$\sum_{m=0}^{\infty} |c_m | 
\big| \tilde{p}_m\big( (\frac{e^{\pi r'}+e^{-\pi r'}}{2})^2\big)\big| $ 
converges.
 From inequality (\ref{ineq:pn}), uniformly absolute convergence of 
 $\sum_{m=0}^{\infty} c_m p_m (\sin ^2 \pi x )$ inside the zone
  $-r'\leq \Im x \leq r'$ is obtained.
Therefore, the proof is complete.
\end{proof}

\subsection{} We note definitions and formulas of elliptic functions.
Let $\omega_1$ and $\omega_3$ be complex numbers such that the value
 $\omega_3/ \omega_1$ is an element of the upper half plane.
The Weierstrass $\wp$-function, the Weierstrass sigma-function and the 
Weierstrass zeta-function are defined as follows:
\begin{equation} \label{wpsigzeta}
\begin{aligned}
\wp (x)&=\wp(x|2\omega_1, 2\omega_3) \\
&= \frac{1}{x^2}+\sum_{(m,n)\in \mathbb{Z} \times \mathbb{Z} \setminus \{ (0,0)\} } 
\Big( \frac{1}{(x-2m\omega_1 -2n\omega_3)^2}
-\frac{1}{(2m\omega_1 +2n\omega_3)^2}\Big), 
\end{aligned}
\end{equation}
\begin{align*}
 \sigma (x)&=x\prod_{(m,n)\in \mathbb{Z} \times \mathbb{Z} \setminus \{(0,0)\} } 
\Big(1-\frac{x}{2m\omega_1 +2n\omega_3}\Big) \\
&\quad\times \exp\Big(\frac{x}{2m\omega_1 +2n\omega_3}+\frac{x^2}{2(2m\omega_1 
+2n\omega_3)^2}\Big), 
\end{align*}
\[
 \zeta(x)=\frac{\sigma'(x)}{\sigma (x)}. 
\]
Setting $\omega_2=-\omega_1-\omega_3$,
$e_i=\wp(\omega_i)$ and  $\eta_i=\zeta(\omega_i)$ for $i=1,2,3$
yields the relations
\begin{gather*}
 e_1+e_2+e_3=\eta_1+\eta_2+\eta_3=0,   \\
 \wp(x)=-\zeta'(x), \quad  (\wp'(x))^2=4(\wp(x)-e_1)(\wp(x)-e_2)(\wp(x)-e_3),  \\
 \wp(x+2\omega_i)=\wp(x), \quad  \zeta(x+2\omega_i)=\zeta(x)+2\eta_i \quad  \; (i=1,2,3), \nonumber \\
 \frac{\wp''(x)}{(\wp'(x))^2}=\frac{1}{2}\Big( \frac{1}{x-e_1}+\frac{1}{x-e_2}
 +\frac{1}{x-e_3} \Big), \\
 \wp(x+\omega_i)=e_i+\frac{(e_i-e_{i'})(e_i-e_{i''})}{\wp(x)-e_i} \quad  
  (i=1,2,3),
\end{gather*}
where $i', i'' \in \{1,2,3\}$ with $i'<i''$, $i\neq i'$ and $i\neq i''$.

The co-sigma functions $\sigma_i(x)$ $(i=1,2,3)$ are 
\begin{equation}
\sigma_i(x)=\exp (-\eta_i x)\sigma(x+\omega_i)/\sigma(\omega _i).
\label{cosigma}
\end{equation}
and satisfy
\begin{equation}
 \big( \frac{\sigma_i(x)}{\sigma(x)}\big)^2 =\wp(x)-e_i, 
\quad(i=1,2,3). \label{sigmai}
\end{equation}
Set $\omega_1=1/2$, $\omega_3=\tau /2$ and $p=\exp (\pi \sqrt{-1} \tau )$. 
The expansion of the Weierstrass $\wp$-function in the variable $p$ is written 
as 
\begin{equation}
\wp (x)=
\frac{\pi^2 }{\sin^2 (\pi x)}- \frac{\pi ^2}{3} -8\pi^2 \sum_{n=1}^{\infty} 
\frac{np ^{2n}}{1-p ^{2n}} (\cos 2n \pi x -1). \label{wpth}
\end{equation}
By setting $x\to  x+1/2$, $x+\tau/2$, $x+(1+\tau)/2$, the following expansions 
are obtained
\begin{equation}\label{wpth1}
\begin{gathered}
 \wp \big(x+\frac{1}{2} \big)=
\frac{\pi^2 }{\cos ^2 (\pi x)}- \frac{\pi ^2}{3} -8\pi^2 \sum_{n=1}^{\infty} 
\frac{np ^{2n}}{1-p ^{2n}} ((-1)^n \cos 2n \pi x -1),  \\
\wp \big(x+\frac{\tau}{2} \big)=- \frac{\pi ^2}{3}
-8\pi^2 \sum_{n=1}^{\infty} np^n\frac{\cos 2\pi n x-p^n}{1-p ^{2n}}, \\
\wp \big(x+\frac{1+\tau}{2} \big)=- \frac{\pi ^2}{3}
-8\pi^2 \sum_{n=1}^{\infty} np^n\frac{(-1)^n\cos 2\pi n x-p^n}{1-p ^{2n}}. 
\end{gathered}
\end{equation}

\subsection*{Acknowledgment}
The author would like to thank Dr. Y. Komori and Prof. T. Oshima for their
fruitful discussions, Prof. A. V. Turbiner, and the anonymous referee for
their valuable comments.
The author was partially supported by the Grant-in-Aid for Scientific Research
(No. 13740021) from the Japan Society for the Promotion of Science.

\end{document}